\documentclass[reqno]{amsart}
\usepackage{graphicx}

\usepackage{amssymb}
\usepackage{amsfonts}
\usepackage{indentfirst}
\usepackage{pst-plot}

\renewcommand\bigskip{\medskip}
\def\@oddhead{\hbox to \textwidth{\footnotesize{\it
Khintchine exponents and Lyapunov exponents } \hfill\thepage}}

\def\bc{\begin{center}}
\def\ec{\end{center}}

\def\a{\alpha}
\def\e{\epsilon}
\def\d{\delta}

\def\r{\theta}

\newtheorem{theorem}{Theorem}[section]
\newtheorem{corollary}[theorem]{Corollary}

\newtheorem{lemma}[theorem]{Lemma}
\newtheorem{proposition}[theorem]{Proposition}

\theoremstyle{definition}

\begin{document}

\title[Generic points  and  Birkhoff ergodic average]
      {Generic points in systems of specification and  Banach valued Birkhoff ergodic average}

\author{Ai-Hua FAN}
\address{Ai-Hua FAN: \
Department of Mathematics,
 Wuhan University,
 Wuhan, 430072, P.R. China
 \& CNRS UMR 6140-LAMFA,
 Universit\'e de Picardie
 80039 Amiens, France}%
\email{ai-hua.fan@u-picaride.fr}%
\author{Ling-Min LIAO}
\address{Ling-Min LIAO: Department of Mathematics,
 Wuhan University,
 Wuhan, 430072, P.R. China \&
 CNRS UMR 6140-LAMFA,
 Universit\'e de Picardie
 80039 Amiens, France}%
 \email{lingmin.liao@u-picardie.fr}%
\author{Jacques PEYRI\`ERE}
\address{Jacques PEYRI\`ERE: Universit\'e Paris-Sud, CNRS UMR
8628, Math\'ematique b\^at.~425, 91405 Orsay Cedex, France}%
 \email{Jacques.Peyriere@math.u-psud.fr}%

\subjclass{Primary: 37B40,37B45; Secondary: 37A05}
 \keywords{Generic points, Systems of specification, Topological entropy, Variational principle}

\begin{abstract}
We prove that systems satisfying the specification property are
saturated in the sense that the topological entropy of the set of
generic points of any invariant measure is equal to the
measure-theoretic entropy of the measure. We study Banach valued
Birkhoff ergodic averages and obtain a variational principle for
 its topological entropy spectrum. As application, we examine a
  particular example concerning with the set of real numbers for
  which the frequencies of occurrences in
their dyadic expansions of infinitely many words are prescribed.
This relies on our explicit determination of a maximal entropy
measure.
\end{abstract}

\maketitle

\section{Introduction}
By dynamical system $(X, T)$, we mean a continuous transformation
$T : X \to X$ on a compact metric space $X$ with metric $d$.  We
shall adopt the notion of topological entropy introduced by Bowen
(\cite{Bowen}, recalled  in the section 2), denoted by $h_{\rm
top}$, to describe the sizes of sets in $X$. We denote by
$\mathcal{M}_{\rm inv}$ the set of all $T$-invariant probability
Borel measures on $X$ and by $\mathcal{M}_{\rm erg}$ its subset of
all ergodic measures. The measure-theoretic entropy of $\mu$ in
$\mathcal{M}_{\rm inv} $ is denoted by $h_\mu$.

Let us first recall some  notions like generic points, saturated
property and the specification property which are quite known in
dynamical systems nowadays.

For $\mu \in \mathcal{M}_{\rm inv}$,  the set $G_{\mu}$ of {\em
$\mu$-generic points}  is defined by
$$
   G_{\mu}:=\left\{x \in X: \frac{1}{n}\sum_{j=0}^{n-1} \delta_{T^j x}
   \stackrel{w^*}{\longrightarrow} \mu \right\},
$$
where $\stackrel{w^*}{\longrightarrow}$ stands for the weak star
convergence of the measures.

A dynamical system $(X, T)$ is said to be {\em saturated} if for
any $\mu \in \mathcal{M}_{\rm inv}$, we have $h_{\rm
top}(G_{\mu})=h_{\mu}$.

Bowen (\cite{Bowen73}) proved that on any dynamical system, we
have $h_{\rm top} (G_\mu)\le h_\mu$ for any $\mu \in
\mathcal{M}_{\rm inv}$. So, saturatedness means that $G_\mu$ is of
optimal topological entropy. One of our main results is to prove
that systems of specification share this saturatedness.

A dynamical system $(X,T)$ is said to satisfy the {\em
specification property} if for any $\e>0$ there exists an integer
$m(\e)\ge 1$ having the property that  for any integer $k\ge2$,
for any $k$ points $x_1,\ldots,x_k$ in $X$, and for any integers
$$
  a_1 \le b_1 < a_2 \le b_2 < \cdots < a_k \le b_k
$$
with $
    a_i - b_{i-1} \ge m(\e) \quad (\forall 2 \le i \le k),
$ there exists a point $y\in X$ such that
$$
   d(T^{a_i+n} y,T^n x_i) < \e
\qquad (\forall
   \ 0 \le n \le b_i-a_i,  \quad \forall 1 \le i \le k).
$$

The specification property was introduced by Bowen
(\cite{Bowen73}) who required that $y$ is periodic. But the
present day tradition doesn't require this. The specification
property implies the topological mixing. Blokh (\cite{Blokh})
proved that these two properties are equivalent for continuous
interval transformations.
Mixing subshifts of finite type satisfy the specification
property. In general, a subshift satisfies the specification if
for any admissible words $u$ and $v$ there exists a word $w$ with
$|w|\le k$ (some constant $k$) such that $uwv$ is admissible. For
$\beta$-shifts defined by $T_\beta x = \beta x (\!\!\mod 1)$,
there is only a countable number of $\beta$'s such that the
$\beta$ shifts admit Markov partition (i.e. subshifts of finite
type), but an uncountable number of $\beta$'s such that the
$\beta$-shifts satisfy the specification property
(\cite{Schmel97}).

Our first result is stated as follows.

\begin{theorem}\label{specification-saturated}
If the dynamical system $(X,T)$ satisfies the specification
property, then it is saturated.
\end{theorem}

As application, we study Banach-valued Birkhoff averages for
saturated systems. Let $\mathbb{B}$ be a real Banach space and
$\mathbb{B}^*$ its dual space, their duality being denoted by
$\langle \cdot, \cdot\rangle$. We consider $\mathbb{B}^*$ as a
locally convex topological space with the weak star topology
$\sigma(\mathbb{B}^*, \mathbb{B})$. For any $\mathbb{B}^*$-valued
continuous function $\Phi: X \to \mathbb{B}^*$, we consider its
Birkhoff ergodic averages
$$
      A_n \Phi(x) = \frac{1}{n} \sum_{j=0}^{n-1} \Phi(T^j x)\qquad (n\ge 1).
$$
We would like to know the asymptotic behavior of $A_n\Phi(x)$ in
the $\sigma(\mathbb{B}^*, \mathbb{B})$-topology for different
points $x \in X$.

Let us state the problem we are studying as follows. Fix a subset
$E \subset \mathbb{B}$. For a sequence $\{\xi_n\} \subset
\mathbb{B}^*$ and a point $\xi \in \mathbb{B}^*$, we denote by
$\limsup_{n \to \infty} \xi_n \stackrel{E}{\le} \xi$ the fact
$$
\limsup_{n \to \infty} \langle \xi_n , w\rangle \le \langle  \xi,
w\rangle\ {\rm \ for \ all\ } w \in E.
$$
The meaning of "$\stackrel{E}{=}$" is obvious. It is clear that
$\limsup_{n \to \infty} \xi_n \stackrel{\mathbb{B}}{\le} \xi$, or
equivalently $\limsup_{n \to \infty} \xi_n
\stackrel{\mathbb{B}}{=} \xi$, means $\xi_n$ converges to $\xi$ in
the weak star topology $\sigma(\mathbb{B}^*, \mathbb{B})$.

Let $\alpha \in \mathbb{B}^*$ and $E \subset \mathbb{B}$. The
object of our study is the set
$$
X_\Phi(\alpha; E) = \left\{ x \in X:  \limsup_{n \to \infty}
A_n\Phi(x) \stackrel{E}{\le}
  \alpha
\right\}.
$$
The set $X_\Phi(\alpha; \mathbb{B})$ will  be simply denoted by
$X_\Phi(\alpha)$. This is the set of points $x \in X$ such that
$\lim_{n \to \infty} A_n\Phi(x) = \alpha$ in $\sigma(\mathbb{B}^*,
\mathbb{B})$-topology. If $\widehat{E}$ denotes the convex cone of
$E$ which consists of all $a w'+b w''$ with $a\ge 0, b\ge 0$ and
$w'\in E, w''\in E$, then $X_\Phi(\alpha, E)=X_\Phi(\alpha,
\widehat{E})$. So we may always assume that $E$ is a convex cone.
If $E$ is symmetric in the sense that $E = -E$, then we have
$$
      X_\Phi(\alpha, E) =\left\{x\in X:
      \lim_{n \to \infty} A_n \Phi(x) \stackrel{E}{=} \alpha\right\}.
$$

 By entropy
spectrum we mean the function
$$
\mathcal{E}^E_\Phi(\alpha) := h_{\rm top} (X_\Phi(\alpha; E)).
$$
Invariant measures will be involved in the study of the entropy
spectrum $\mathcal{E}^E_\Phi(\alpha)$.  We set
$$
       \mathcal{M}_\Phi(\alpha; E) = \left\{\mu\in \mathcal{M}_{\rm
       inv}: \int \Phi d\mu \stackrel{E}{\le} \alpha \right\}
$$
where $\int \Phi d \mu$ denotes the vector-valued integral in
Pettis' sense (see \cite{Rudin}) and the inequality
"$\stackrel{E}{\le}$" means
$$
             \int \langle \Phi, w \rangle d\mu \le
             \langle \alpha, w \rangle
             \ \ \mbox{\rm for \ all }\ w \in E.
$$

For saturated systems, we prove the following variational
principle.

\begin{theorem}\label{VP}
Suppose that the dynamical system $(X, T)$ is saturated. Then
\\
\indent {\rm (a) }\ If  $
\mathcal{M}_\Phi(\alpha; E)= \emptyset$, we have $X_\Phi(\alpha, E)=\emptyset$.\\
\indent {\rm (b)}\ If  $ \mathcal{M}_\Phi(\alpha; E)\not=
\emptyset$, we have
\begin{equation}\label{variational-principle}
         h_{\rm top} (X_\Phi(\alpha; E)) =
        \sup_{\mu \in \mathcal{M}_\Phi(\alpha; E)} h_\mu.
\end{equation}
\end{theorem}

When $\mathbb{B}$ is a finite dimensional Euclidean space
$\mathbb{R}^d$ and $E= \mathbb{R}^d$, the variational principle
(\ref{variational-principle}) with $E= \mathbb{R}^d$ was proved in
\cite{FF,FFW} for subshifts of finite type, then for conformal
repellers (\cite{FLW}) and later generalized to systems with
specification property \cite{TV}.  There are other works assuming
that $\Phi$ is regular (H\"{o}lder for example). See
\cite{Besicovitch,Eggleston} for classical discussions,
\cite{BSS,BSS2, Oli1,Oli2, PS, Temp} for recent developments on
Birkhoff averages, and \cite{BBP, BBH, BMP, CLP,F} for the
multifractal analysis of measures.

The study of infinite dimensional Birkhoff averages is a new
subject. We point out that \cite{Pey} provides  another point of
view, i.e. the thermodynamical point of view which was first
introduced by physicists.

The above variational principle~(\ref{variational-principle}),
when $E=\mathbb{B}$, is easy to generalize to the following
setting. Let $\Psi$ be a continuous function defined on the closed
convex hull of the image $\Phi(X)$ of $\Phi$ into a topological
space~$Y$. For given $\Phi$, $\Psi$, and $\beta\in Y$, we set
\begin{equation*}
X_\Phi^\Psi(\beta) = \left\{ x\in X\ :\ \lim_{n\to\infty}
\Psi\bigl( A_n\Phi(x)\bigr) = \beta \right\}
\end{equation*}
and
\begin{equation*}
{\mathcal M}_\Phi^\Psi(\beta) = \left\{ \mu\in {\mathcal M}_{\rm
inv}\ :\ \Psi\left( \int \Phi\,d\mu \right) = \beta\right\}.
\end{equation*}
We also set
\begin{equation*}
\widehat{X}_\Phi^\Psi(\beta) = \left\{ x\in X\ :\ \Psi\left(
\lim_{n\to\infty} A_n\Phi(x)\right) = \beta  \right\} =
\bigcup_{\alpha:\,\Psi(\alpha)=\beta} X_\Phi(\alpha).
\end{equation*}
It is clear that $\widehat{X}_\Phi^\Psi(\beta)$ is a subset
of~${X}_\Phi^\Psi(\beta)$.

\begin{theorem}\label{saturated}
Suppose that the dynamical system $(X,T)$ is saturated. Then
\begin{enumerate}
\item if ${\mathcal M}_\Phi^\Psi(\beta)=\emptyset$, we have
$X_\Phi^\Psi(\beta)=\emptyset$, \item if ${\mathcal
M}_\Phi^\Psi(\beta)\ne\emptyset$, we have
\begin{equation}\label{g-principle}
h_{\rm top}\bigl( X_\Phi^\Psi(\beta)\bigr) = h_{\rm top}\bigl(
\widehat{X}_\Phi^\Psi(\beta)\bigr) = \sup_{\mu\in {\mathcal
M}_\Phi^\Psi(\beta)} h_\mu.
\end{equation}
\end{enumerate}
\end{theorem}

This generalized variational principle~(\ref{g-principle}) allows
us to study generalized ergodic limits like
\begin{equation}
\lim_{n\to\infty}\frac{\sum_{j=0}^{n-1}\Phi(T^jx)}{\sum_{j=0}^{n-1}g(T^jx)
},
\end{equation}
where $g: X \to \mathbb{R}^+$ is a continuous positive function.
It suffices to apply~(\ref{g-principle}) to~$\Phi$ replaced by
$(\Phi,g)$ and~$\Psi$ defined by $\Psi(x,y) = x/y$, with
$x\in{\mathbb B}^*$ and $y\in{\mathbb R}^+$.

It also allows us to study the set of points $x\in X$ for which
the limits $A_\infty\Phi(x) = \lim_{n\to\infty}A_n\Phi(x)$ verify
the equation
\begin{equation}
\Psi\bigl( A_\infty\Phi(x)\bigr) = \beta.
\end{equation}


The present infinite dimensional version of variational principle
would have many interesting applications. We will just illustrate
the usefulness of the variational principle by the following study
of frequencies of blocks in the dyadic development of real
numbers. It can be reviewed as an infinitely multi-recurrence
problem.

Let us state the question to which we can answer. All but a
countable number of real numbers $t \in [0, 1]$ can be uniquely
developed as follows
$$
     t = \sum_{n=1}^\infty \frac{t_n}{2^n}
     \qquad (t_n \in \{0, 1\}).
$$
Let $k\ge 1$. We write $0^k$ for the block of $k$ consecutive
zeroes and we define the $0^n$-frequency of $t$ as the limit (if
it exists)
$$
     f(t,k) = \lim_{n\to\infty} \frac{1}{n}\sharp\{1\le j \le n:
          t_j t_{j+1}\cdots t_{j+k-1} = 0^k\}.
$$
Let $(a_1, a_2, \cdots)$ be a sequence of non-negative numbers. We
denote by $S(a_1, a_2, \cdots)$ the set of all numbers $t\in
[0,1]$ such that $f(t, k) = a_k$ for all $k\ge 1$. As a
consequence of the variational principle
(~\ref{variational-principle}), we prove

\begin{theorem} The set
$S(a_1, a_2, \cdots)$ is non-empty if and only if the following
condition is satisfied
\begin{equation}\label{convex-condition}
   1=a_0 \ge a_1 \ge a_2 \ge \dots; \quad a_i - 2 a_{i+1} +
    a_{i+2}\ge 0  \ \ \ \  (i\ge 0).
\end{equation}
If the condition~(\ref{convex-condition}) is fulfilled, we have
\begin{equation}
  h_{\rm top}
  (S(a_1, a_2, \cdots)) = -h(1-a_1) + \sum_{j=0}^\infty h(a_j - 2a_{j+1}+a_{j+2})
\end{equation}
where $h(x) = -x \log x$.
\end{theorem}

Furthermore, it is proved that there is a unique maximal entropy
measure, which is completely determined (see Lemma~\ref{unicite}).

The paper is organized as follows. In the section~2, we give some
preliminaries. In the section~3, we prove Theorem
\ref{specification-saturated}. In the section~4, we prove the
theorems~\ref{VP} and~\ref{saturated} and examine the case where
${\mathbb B}=\ell^1({\mathbb Z})$. In the section~5, we apply the
variational principle~(\ref{variational-principle}) to the study
of the recurrence into an infinite number of cylinders of the
symbolic dynamics. Especially, we study the set of orbits whose
recurrences into infinitely many cylinders are prescribed. This
relies on the explicit determination of a maximal entropy measure
which, by definition, maximizes the supremum
in~(\ref{variational-principle}).

\section{Preliminary}
Before proving the main theorems, we wish to recall the notions of
topological entropy, the Bowen lemma and two propositions about
the measure-theoretic entropy.

Recall that $X$ is a compact metric space with its metric $d$ and
that $T: X \to X$ is a continuous transformation on $X$. For any
integer $n\ge 1$ we define the Bowen metric $d_n$ on $X$ by
$$
      d_n(x, y) = \max_{0\le j <n} d(T^jx, T^j y).
$$
For any $\epsilon >0$, we will denote by $B_n(x, \epsilon)$ the
open $d_n$-ball centered at $x$ of radius $\epsilon$.

\subsection{Topological entropy and Bowen lemma}\

Let $ Z \subset X$ be a subset of $X$. Let $\epsilon >0$. We say
that a collection (at most countable) $R=\{B_{n_i}(x_i,
\epsilon)\}$ covers $Z$ if $Z \subset \bigcup_i B_{n_i}(x_i,
\epsilon)$.  For such a collection $R$, we put $n(R) = \min_i
n_i$.  Let $s\ge 0$. Define
$$
     H^s_n (Z, \epsilon) = \inf_R \sum_i \exp(-s n_i),
$$
where the infimum is taken over all covers $R$ of $Z$ with $n(R)
\ge n$.  The quantity $H^s_n (Z, \epsilon)$ is non-decreasing as a
function of $n$, so the following limit exists
$$
   H^s (Z, \epsilon)  = \lim_{n \to \infty} H^s_n (Z, \epsilon).
$$
For the quantity $H^s (Z, \epsilon)$ considered as a function of
$s$, there exists a critical value, which we denote by $h_{\rm
top} (Z, \epsilon)$, such that
$$
     H^s (Z, \epsilon) =\left\{ \begin{array} {ll} +\infty, & s <
                        h_{\rm top} (Z, \epsilon) \\ 0 , & s> h_{\rm
                        top} (Z, \epsilon).  \end{array} \right.
$$
One can prove that the following limit exists
$$
          h_{\rm top} (Z) = \lim_{\epsilon \to 0} h_{\rm top} (Z,
          \epsilon).
$$
The quantity $h_{\rm top} (Z)$ is called the {\em topological
entropy} of $Z$ (\cite{Bowen}).

For $x \in X$, we denote by $V(x)$ the set of all weak limits of
the sequence of probability measures $n^{-1}\sum_{j=0}^{n-1}
\delta_{T^j x}$. It is clear that $V(x) \not=\emptyset$ and $V(x)
\subset \mathcal{M}_{\rm inv}$ for any $x$. The following Bowen
lemma is one of the key lemmas for proving the variational
principle.

\begin{lemma}[Bowen \cite{Bowen73}]\label{Bowen}
    For $t\ge 0$, consider the set
    $$
         B^{(t)} =\left\{x\in X:  \exists \mu \in V(x) \
         \mbox{\rm satisfying}\ h_\mu \le t\right\}.
    $$
    Then $h_{\rm top} (B^{(t)})\le t$.
\end{lemma}

Let $\mu \in \mathcal{M}_{\rm inv}$ be an invariant measure. A
point $x \in X$ such that $V(x)=\{\mu\}$ is said to be generic for
$\mu$. Recall our definition of $G_\mu$, we know that $G_\mu$ is
the set of all generic points for $\mu$. Bowen proved that $h_{\rm
top}(G_\mu)\le h_\mu$ for any invariant measure. This assertion
can be deduced by using Lemma~\ref{Bowen}. In fact, the reason is
that $x \in G_\mu $ implies $\mu \in V(x)$. Bowen also proved that
the inequality becomes equality when $\mu$ is ergodic. However, in
general, we do not have the equality $ h_{\rm top} (G_\mu) =
h_\mu$ (saturatedness) and it is even possible that $G_\mu =
\emptyset$. Cajar \cite{Cajar} proved that full symbolic spaces
are saturated. Concerning the $\mu$-measure of $G_\mu$, it is well
known that $\mu(G_\mu) = 1 \ {\rm or}\ 0$ according to whether
$\mu$ is ergodic or not (see \cite{DGS}).

 \setcounter{equation}{0}

\subsection{Two propositions about the measure-theoretic entropy}\

 We denote by $C(X)$  the set of continuous
functions on $X$, by $\mathcal{M}= \mathcal{M}(X)$ the set of all
Borel probability measures.

In the sequel, we fix a sequence $(p_i)_{i \ge 1}$ such that $p_i
> 0$ for all $i \ge 1$ and $\sum_{i=1}^\infty p_i =1$ (for example,
$p_i=2^{-i}$ will do). Suppose that $s_n = (s_{n, i})_{i \ge 1}$
($n=1,2, \cdots$) is a sequence of elements in $\ell^\infty$.
 It is obvious that $s_n$
converges to $\alpha=(\alpha_i)_{i \ge 1} \in \ell^{\infty}$ in
the weak star topology (i.e. each coordinate converges) is
equivalent to
$$
      \lim_{n \to \infty} \sum_{i=1}^\infty p_i |s_{n,
      i}-\alpha_i|=0.
$$
We also fix a sequence  of continuous functions
$\{\Phi_1,\Phi_2,\ldots\}$ which is dense in the unit ball of
$C(X)$. Write $\Phi=(\Phi_1,\Phi_2,\ldots)$. It is evident that
$\Phi: X \to \ell^{\infty}$ is continuous when $\ell^{\infty}$ is
equipped with its weak star topology. Fix an invariant measure
$\mu \in \mathcal{M}_{\rm inv}$. Let
$$
     \alpha = (\alpha_1, \alpha_2,  \cdots \alpha_i, \cdots )
     \quad {\rm where} \ \  \alpha_i = \int \Phi_i d \mu.
$$
The set of generic points $G_\mu$ can be described as follows
\begin{eqnarray}
G_\mu = \left \{ x\in X:\lim_{n \to \infty} \sum_{i=1}^{\infty}
p_i \left| A_n \Phi_i - \alpha_i \right|=0 \right\} =
X_\Phi(\alpha).
\end{eqnarray}

It is well known that the weak topology of $\mathcal{M}$ is
compatible with the topology induced by the metric
\begin{eqnarray}
 \tilde{d}(\mu, \nu)= \sum_{i=1}^\infty p_i \left|\int \Phi_i d \mu
               -\int\Phi_i d \nu \right|
\end{eqnarray}
where both $(p_i)_{i \ge 1}$ and $\{\Phi_i\}_{i\ge 1}$ are chosen
as above.

The following two results will be useful for us.

\begin{proposition}[Young  \cite{Young}]\label{Young}
For any  $\mu \in \mathcal{M}_{\rm inv}$ and any
 numbers $0<\delta<1$ and  $0<\theta<1$, there exist an
invariant measure $\nu$ which is a finite convex combination  of
ergodic measures, i.e.
$$\nu =\sum\limits_{k=1}^{r}\lambda_{k}\nu_{k},
\qquad {\rm where} \ \ \lambda_{k}>0,\ \
\sum\limits_{k=1}^{r}\lambda_{k}=1,\ \  \nu_{k}\in
\mathcal{M}_{\rm erg}, \ r\in \mathbb{N}^+
$$
such that $$ \tilde{d}(\mu, \nu) <\delta,\qquad h_\nu\geq
h_\mu-\theta.
$$
\end{proposition}

This is a consequence of the following result due to Jacobs (see
\cite{Walters}, p. 186). Let $\mu \in \mathcal{M}_{\rm inv}$ be an
invariant measure which has the ergodic decomposition $\mu =
\int_{\mathcal{M}_{\rm erg}} \tau d \pi(\tau)$ where $\pi$ is a
Borel probability measure on $\mathcal{M}_{\rm erg}$. Then we have
$$
   h_{\mu}= \int_{\mathcal{M}_{\rm erg}} h_{\tau} d \pi(\tau).
$$

\begin{proposition} [Katok  \cite{Katok}]\label{Katok}
 Let $\mu \in \mathcal{M}_{\rm erg}$ be an ergodic invariant measure. For $\e>0$ and $\d>0$,
let $r_n(\e,\d,\mu)$ denote the minimum number of $\e$-balls in
the Bowen metric $d_n$  whose union has $\mu$-measure more than or
equal to $1-\d$. Then for each $\d>0$ we have
$$
    h_\mu=\lim_{\e\to 0}\limsup_{n\to \infty}\frac{1}{n}\log r_n(\e,\d,\mu)=
    \lim_{\e\to 0}\liminf_{n\to \infty}\frac{1}{n}\log
    r_n(\e,\d,\mu).
$$
\end{proposition}

 In \cite{Katok}, it was assumed that $T:X\rightarrow X$ is a homeomorphism.
 But the proof in \cite{Katok} works for the transformations we are studying.

\section{Systems with specification property are saturated}
\setcounter{equation}{0}
 In this section, we prove Theorem \ref{specification-saturated} which
says that every system satisfying the  specification property is
saturated. Because of Bowen's lemma (Lemma~\ref{Bowen}), we have
only to show $h_{\rm top} (G_\mu) \ge h_\mu$.  The idea of the
proof appeared in \cite{FF, FFW} and was developed in \cite{TV}.
It consists of constructing the so-called dynamical Moran sets
which approximate the set of generic points $G_\mu$.

\subsection{Dynamical Moran sets and their entropies}\

Fix $\epsilon >0$.
 Let $\{m_k\}_{k\ge 1}$ be the sequence of integers defined by $m_k=m(2^{-k}\e)$
 which is the constant appeared in the definition
 of the
specification property ($k=1,2,\ldots$). Let $\{W_k\}_{k\ge 1}$
 be a sequence of finite sets in $X$ and
$\{n_k\}_{k\ge 1}$
 be a sequence of positive integers. Assume that
\begin{equation}\label{2}
   d_{n_k}(x,y) \geq 5\e \qquad (\forall x,y \in W_k \quad x\neq
   y).
\end{equation}
Let $\{N_k\}_{k \geq 1}$ be another sequence of positive integers
with $N_1=1$. Using these data, we are going to construct a
compact set of Cantor type, called a dynamical Moran set, which
will be denoted by $F= F(\epsilon, \{W_k\}, \{n_k\}, \{N_k\})$. We
will give an estimate for its topological entropy.

Denote
$$M_k= \# W_k.$$
Fix $k\ge 1$.  For any $N_k$ points $x_{1},\cdots,x_{N_k}$ in
$W_k$ i.e. $(x_{1},\cdots,x_{N_k})\in W_k^{N_k}$, we choose a
point $y(x_{1},\cdots,x_{N_k}) \in X$, which does exist by the
specification property, such that
\begin{equation}\label{2bis}
   d_{n_k}(x_{s},T^{a_s}y) < \frac{\e}{2^k} \qquad (s=1,\ldots,
   N_k)
\end{equation} where
$$
   a_s=(s-1)(n_k + m_k).
$$
Both (\ref{2}) and (\ref{2bis}) imply that for two distinct points
$(x_{1},\cdots,x_{{N_k}})$ and $
(\bar{x}_{1},\cdots,\bar{x}_{{N_k}})$ in $W_k^{N_k}$ we have
\begin{equation}\label{2ter}
   d_{t_k}(y(x_{1},\cdots,x_{{N_k}}),
   y(\bar{x}_{1},\cdots,\bar{x}_{{N_k}})) > 4\e
\end{equation} where $t_k = a_{N_k} + n_k$, i.e.
$$
   t_k=  (N_k-1)m_k + N_k n_k.
$$
In fact, let $y = y(x_1, \cdots, x_{N_k})$ and $\bar{y}=
y(\bar{x}_{1},\cdots,\bar{x}_{{N_k}})$. Suppose $x_{s} \neq
\bar{x}_{s}$ for some $s\in\{1,\cdots,N_k\}$. Then
\begin{eqnarray*}
d_{t_k}(y,\bar{y}) & \ge &
   d_{n_k}(T^{a_s}y,T^{a_s}\bar{y})\\
  & \geq & d_{n_k}(x_{s}, \bar{x}_{s}) -
   d_{n_k}(x_{s},T^{a_s}y)
        -d_{n_k}(\bar{x}_{s},T^{a_s}\bar{y})\\
   &>& 5\e-\e/2- \e/2 = 4\e.
\end{eqnarray*}

Let
$$ D_1=W_1, \quad
D_k=\left\{ y(x_{1},\cdots,x_{N_k}): (x_{1},\cdots,x_{{N_k}}) \in
W_k^{N_k}  \right\} \ \ (\forall k \ge 2).
$$
Now define recursively $L_k$ and $\ell_k$ as follows. Let
$$
   L_1=D_1, \quad \ell_1=n_1.
$$
 For any $x\in L_k$ and any $y\in
D_{k+1} \ (k\geq 1)$, by the specification property, we can find a
point $z(x,y) \in X$ such that
$$
    d_{\ell_k}(z(x,y),x)<\frac{\e}{2^{k+1}}, \qquad
    d_{t_{k+1}}(T^{\ell_k+m_{k+1}}z(x,y),y)<\frac{\e}{2^{k+1}}.
$$
We will choose one and only one such $z(x,y)$ and call it the
descend from $x\in L_k$ through $y \in D_{k+1}$. Let
\begin{eqnarray*}
       L_{k+1}
       &=&\left\{z(x, y):   x \in L_k, y \in
       D_{k+1}\right\},\\
\ell_{k+1}& =& \ell_k + m_{k+1} + t_{k+1}=N_1n_1 +\sum_{i=2}^{k+1}
N_i(m_i+n_i).
\end{eqnarray*}
Observe that for any $x \in L_k$ and for all $y,\bar{y}\in
D_{k+1}$ with $ y\neq \bar{y}$, we have
\begin{equation}\label{Lk1}
   d_{\ell_k}(z(x,y),z(x,\bar{y})) < \frac{\e}{2^k} \quad (k\geq 1),
\end{equation}
and for any $x, \bar{x} \in L_k $ and $y, \bar{y} \in D_{k+1}$
with $(x,y)\neq (\bar{x},\bar{y})$, we have
\begin{equation}\label{Lk2}
     d_{\ell_{k+1}}(z(x,y),z(\bar{x},\bar{y})) > 3\e \quad (k\geq 1).
\end{equation}
The fact (\ref{Lk1}) is obvious. To prove (\ref{Lk2}), first
remark that $d_{\ell_1}(z,\bar{z})\ge 5\e > 4\e$ for any $z,
\bar{z}\in L_1$ with $z \neq \bar{z}$, and that for any $x,\bar{x}
\in L_k$ and $y,\bar{y} \in D_{k+1}$ with $(x,y) \neq
(\bar{x},\bar{y})$ we have
$$
    d_{\ell_{k+1}}(z(x,y),z(\bar{x},\bar{y}))
    \geq  d_{\ell_k}(x,\bar{x}) - d_{\ell_k}(z(x,y),x)
    -d_{\ell_k}(z(\bar{x},\bar{y}),\bar{x})
$$
and
\begin{eqnarray*}
  & & d_{\ell_{k+1}}(z(x,y),z(\bar{x},\bar{y}))\\
   &\ge&
   d_{t_{k+1}}(y,\bar{y}) - d_{t_{k+1}}(T^{\ell_k + m_{k+1}}z(x,y),y)
    -d_{t_{k+1}}(T^{\ell_k + m_{k+1}}z(\bar{x},\bar{y}),\bar{y}).
\end{eqnarray*}
Now using the above two inequalities, we prove (\ref{Lk2}) by
induction.
 For any $x,\bar{x} \in
L_1$ and $y,\bar{y} \in D_{2}$ with either $x\neq \bar{x}$ or
$y\neq \bar{y}$, we have
$$
   d_{\ell_2}(z(x,y),z(\bar{x},\bar{y})) > 4\e- \frac{\e}{2^2}-
   \frac{\e}{2^2} = 4\e- \frac{\e}{2}.
$$
Suppose we have obtained that
$$
    d_{\ell_{k}}(z(x,y),z(\bar{x},\bar{y})) > 4\e-
    \frac{\e}{2}- \frac{\e}{2^2}- \cdots - \frac{\e}{2^{k-1}}.
$$
Then for any $x,\bar{x} \in L_k$ and $y,\bar{y} \in D_{k+1}$ with
$(x,y) \neq (\bar{x},\bar{y})$ we have
\begin{eqnarray*}
    d_{\ell_{k+1}}(z(x,y),z(\bar{x},\bar{y}))
    &>& 4\e-\frac{\e}{2}- \frac{\e}{2^2}- \cdots - \frac{\e}{2^{k-1}} -
        \frac{\e}{2^{k+1}} - \frac{\e}{2^{k+1}}\\
    &=& 4\e-\frac{\e}{2}- \frac{\e}{2^2}- \cdots - \frac{\e}{2^{k-1}} -
        \frac{\e}{2^{k}}
    > 3\e
\end{eqnarray*}

Now define our dynamical Moran set
$$
     F=F(\epsilon, \{W_k\},\{n_k\},\{N_k\})=\bigcap_{k=1}^\infty F_k,
     $$
     where
     $$ F_k=\bigcup_{x \in L_k} \overline{B}_{\ell_k}(x,\e2^{-(k-1)})
$$
($\overline{B}(x, r)$ denoting the closed ball of center $x$ and
radius $r$). The set $F$ is Cantor-like because for any distinct
points $ x', x'' \in L_k $, by (\ref{Lk2}) we have
$$
    \overline{B}_{\ell_k}(x', \e2^{-(k-1)}) \bigcap
    \overline{B}_{\ell_k}(x'', \e2^{-(k-1)}) = \emptyset
$$
and if $z \in L_{k+1}$ descends from $x \in L_k$, by (\ref{Lk1})
we have
$$
     \overline{B}_{\ell_{k+1}}(z, \e2^{-k}) \subseteq
     \overline{B}_{\ell_k}(x,\e2^{-(k-1)}).
$$

\begin{proposition} [Entropy of $F$]\label{entofF} For any integer
$n\ge 1$,  let $k=k(n)\ge 1 $ and $0\le p=p(n)< N_{k+1}$ be the
unique integers such that
$$
   \ell_k + p(m_{k+1}+ n_{k+1}) < n \leq \ell_k + (p+1)(m_{k+1}+
   n_{k+1}).
$$
We have
$$
h_{\rm top}(F) \geq \liminf_{n \to \infty} \frac{1}{n}(N_1 \log
M_1 +\cdots +N_k \log M_k + p \log M_{k+1}).
$$
\end{proposition}

 \proof
For every $k\ge 1$, consider the discrete measure $\sigma_k$
concentrated on $F_k$
$$
      \sigma_k = \frac{1}{\# L_k} \sum_{ x \in L_k} \delta_x.
$$
It can be proved that
  $\sigma_k$ converges in the week star
  topology to a probability measure $\sigma$ concentrated on $F$.
  Moreover, for sufficiently large $n$ and every point $x\in X$
such that $
         B_n(x,\e/2) \cap F \neq \emptyset,
$ we have
$$
         \sigma(B_n(x,\e/2)) \leq \frac{1}{\#(L_k)M_{k+1}^p} =
         \frac{1}{M_1^{N_1} \cdots M_k^{N_k} M_{k+1}^p}.
$$
 (see \cite{TV}). Then we apply the mass distribution principle to estimate the
 entropy. $\Box$

\subsection{Box-counting of $G_\mu$}\

Recall that $\alpha=(\alpha_i)_{i \ge 1} \in \ell^{\infty}$ and
$\Phi=(\Phi_i)$ is a dense sequence in the unit ball of $C(X)$.
For $\delta >0$ and $n \ge 1$, define
$$
X_\Phi(\alpha, \delta, n) =\left\{x \in X:
      \sum_{i=1}^\infty p_i |A_n \Phi_i(x)
      -\alpha_i| <\delta
      \right\}.
$$
For $\epsilon>0$, let $N(\alpha,\d,n,\e)$ denote the minimal
number of balls $B_{n}(x,\e)$ to cover the set
$X_\Phi(\alpha,\delta,n)$. Define
\begin{equation}\label{Lambda}
\Lambda_\Phi (\alpha):=\lim\limits_{\e\to 0}\lim\limits_{\d\to
0}\limsup\limits_{n\to \infty}\frac{1}{n} \log N(\alpha,\d,n,\e)
\end{equation}
By the same argument in \cite{FF} (p. 884-885), we can prove the
existence of the limits, and the following equality:
$$
\Lambda_\Phi (\alpha)=\lim\limits_{\e\to 0}\lim\limits_{\d\to
0}\liminf\limits_{n\to \infty}\frac{1}{n} \log N(\alpha,\d,n,\e).
$$

\begin{proposition}\label{lamb>hu}
  $\Lambda_\Phi (\alpha)\geq h_\mu$.
\end{proposition}

\proof It suffices to prove $\Lambda_\Phi (\alpha) \ge h_\mu -
4\theta$ for any $\theta > 0$. For each $i\geq 1$, define the
variation of $\Phi_i$ by
$$
     {\rm var}(\Phi_i, \epsilon) = \sup_{d(x, y)<\epsilon} |\Phi_i(x)-
     \Phi_i(y)|.
$$
By the compactness of $X$ and the continuity of $\Phi_i$,
$\lim_{\epsilon \to 0}{\rm var}(\Phi_i,\e)\rightarrow 0$. So
$$ \lim_{\epsilon \to 0}
\sum_{i=1}^{+\infty} p_i {\rm var}(\Phi_i,\e)\rightarrow 0.
$$
This, together with (\ref{Lambda}), allows us to choose $\e>0$ and
$\d>0$ such that
\begin{equation}\label{3}
  \sum_{i=1}^{+\infty} p_i {\rm var}(\Phi_i,\e)<\d<\theta
\end{equation}
and
\begin{equation}\label{4}
 \limsup_{n\to \infty}\frac{1}{n}\log N(\a,5\d,n,\e) < \Lambda_\Phi
 (\alpha)+\theta.
\end{equation}

For the measure $\mu$, take an invariant measure $\nu=\sum_{k=1}^r
\lambda_k \nu_k$ having the properties stated  in Proposition
\ref{Young}. For $1\le k\le r$ and $N\ge 1$, set
$$
   Y_k(N)= \left\{x\in X:
            \sum_{i=1}^\infty p_i \left|A_n\Phi_i(x) - \int \Phi_i d
             \nu_k\right|<\delta \ \ \  (\forall n \ge N)\right\}.
$$
Since $\nu_k$ is ergodic,  by the Birkhoff theorem, we have
\begin{equation}\label{5}
     \lim_{n\to \infty} \sum_{i=0}^{\infty} p_i \left| A_n\Phi_i(x) - \int \Phi_i d
     \nu_k \right| = 0
     \quad \nu_k{\rm -a.e.}
\end{equation}
Then by the Egorov theorem, there exists a set with
$\nu_k$-measure greater than $1-\r$ on which the  above limit
(\ref{5}) is uniform. Therefore, if $N$ is sufficiently large, we
have
\begin{equation}\label{5bis}
        \nu_k(Y_k(N))>1-\r   \qquad (\forall\  k=1,\cdots,r).
\end{equation}

Apply the second equality in Proposition \ref{Katok} to the triple
$(\nu_k, 4\epsilon, \theta)$ in place of $(\mu, \epsilon,
\delta)$. When $\epsilon >0 $ is small enough, we can find an
integer $N_k =N_k(\nu_k, 4\epsilon, \theta)\ge 1$ such that
$$
    r_{n}(4\e,\r,\nu_k) \geq \exp (n (h_{\nu_k} - \r))
     \qquad (\forall n \ge N_k).
$$
This implies that if $n \ge N_k$, then the minimal number of balls
$B_{n}(x, 4\e)$ to cover $Y_{k}(N)$ is greater than or equal to
$\exp(n (h_{\nu_{k}}-\r))$. Consequently, if we use $C(n,4\e)$ to
denote a maximal $(n,4\e)$-separated set in $Y_{k}(N)$, then
\begin{equation}\label{6}
\#C(n,4\e)\geq\exp(n (h_{\nu_{k}}-\r))
 \qquad (\forall n \ge N_k).
\end{equation}

 Choose a sufficiently large integer
$N_{0}$ such that
$$
     n_{k}:=[\lambda_{k}n] \ge \max(N_1, \cdots,N_{k},N)
    \qquad (\forall
k=1,\ldots,r; \ \forall n\ge N_{0})
$$
($[\cdot]$ denoting the integral part). By the specification
property, for each  $r$ points $x_1 \in C(n_1 , 4\e),\ldots,x_r\in
C(n_r , 4\e)$, there exist an integer $m(\e)$ depending on $\e$
and a point $y=y(x_1,\ldots,x_r) \in X$ such that
\begin{equation}\label{7}
d_{n_k}(T^{a_k}y, x_k)<\e \qquad (1\le k\le r)
\end{equation}
where
$$
a_1=0,\ \  a_k =(k-1)m+\sum\limits_{s=1}^{k-1}n_s \ \ (k\geq 2).
$$
Write $\hat{n}=a_r + n_r$, i.e.
$$
\hat{n}=(r-1)m+ \sum_{s=1}^{r}n_s.
$$
We claim that   for all such $y=y(x_1,\ldots,x_r)$, we have
\begin{equation}\label{claim1}
     y=y(x_1,\ldots,x_r) \in X_\Phi(\alpha, 5 \delta, \hat{n})
\end{equation}
when $n$ is sufficiently large, and that for two distinct points
$(x_1,\ldots,x_r)$ and $(x_1',\ldots,x_r')$ in $C(n_1, 4 \epsilon)
\times \cdots \times C(n_r, 4 \e)$, the points
$y=y(x_1,\ldots,x_r)$ and $y'=y(x_1',\ldots,x_r')$ satisfy
\begin{equation}\label{claim2}
    d_{\hat{n}}(y,y') > 2\e.
\end{equation}

If we admit (\ref{claim1}) and (\ref{claim2}), we can conclude. In
fact, the balls $B_{\hat{n}}(y,\e)$ are disjoint owing to
(\ref{claim2}) and hence there are $\#C(n_1,4\e) \times \ldots
\times \#C(n_r,4\e)$ such balls. Therefore, because of
(\ref{claim1}), the minimal number of ($\hat{n},\e)$-balls needed
to cover $X_\Phi(\alpha,5\d,\hat{n})$ is greater than the number
of such points $y$'s. That is to say
$$
  N(\a,5\d,\hat{n},\e) \geq \#C(n_1,4\e) \times \ldots \times \#C(n_r,4\e)
$$
Then by (\ref{6}), we get
$$
   N(\a,5\d,\hat{n},\e) \geq \exp \sum_{k=1}^r [\lambda_k
   n](h_{\nu_k}-\r).
$$
By noticing that $\frac{[\lambda_k
   n]}{\hat{n}} \to \lambda_k$ as $n \to \infty$ and $\sum_{k=1}^r
   \lambda_k
   =1$, we get
$$
   \liminf_{\hat{n}\to \infty}\frac{1}{\hat{n}} \log N(\a,5\d,\hat{n},\e) \geq h_\mu-3\r
$$
This, together with (\ref{4}), implies $
  \Lambda_\Phi(\a)\geq h_\mu-4\r
$.

Now return to prove (\ref{claim1}) and (\ref{claim2}). The proof
of (\ref{claim2}) is simple: suppose $x_k \neq x_k'$ for some
$1\le k\le r$. By (\ref{7}),
$$
    d_{\hat{n}}(y,y')\geq d_{n_k}(T^{a_k}y,T^{a_k}y') \geq d_{n_k}(x_k,x_k')-2\e > 4\e-2\e=
    2\e.
$$
Now prove (\ref{claim1}). Recall that $\alpha_i = \int \Phi_i d
\nu$ and $\nu = \sum_{k=1}^r \lambda_k \nu_k$. We have
$$
    \left|A_{\hat{n}}\Phi_i(y) - \a_i \right|
    \le
     \left|A_{\hat{n}}\Phi_i(y) - \sum_{k=1}^r \lambda_k
                       \int\Phi_{i} d\nu_k  \right| +
                                \left| \int\Phi_{i} d\nu - \int\Phi_{i} d\mu
                                \right|.
$$
Since $\tilde{d}(\mu, \nu)<\delta$ i.e. $\sum_{i=1}^\infty p_i
|\int \Phi_i d \mu - \int \Phi_i d \nu|<\delta$, we have only to
show that
\begin{equation}\label{last}
   \sum_{i=1}^\infty p_i \left|A_{\hat{n}}\Phi_i(y) -
                      \sum_{k=1}^r \lambda_k  \int\Phi_{i} d\nu_k  \right| <4 \delta.
\end{equation}
Write
\begin{eqnarray*}
    A_{\hat{n}} \Phi_i(y)
      & = & \frac{1}{\hat{n}}\sum_{k=1}^r \sum_{j=0}^{[\lambda_k n]-1}
      \Phi_i(T^{a_k+j}y)
         + \frac{1}{\hat{n}} \sum_{k=2}^r \sum_{j=a_k-m}^{a_k -1} \Phi_i(T^j
         y)\\
             & = & \sum_{k=1}^r  \frac{[\lambda_k n]}{\hat{n}} A_{[\lambda_k
             n]}\Phi_i(T^{a_k}y) + \frac{1}{\hat{n}} \sum_{k=2}^r \sum_{j=a_k-m}^{a_k -1} \Phi_i(T^j
         y).
         \end{eqnarray*}
Then
\begin{eqnarray*}
   \left|A_{\hat{n}}\Phi_i(y) - \sum_{k=1}^r \lambda_k
                       \int\Phi_{i} d\nu_k  \right| \le I_1(i) +
                       I_2(i) + I_3(i) + I_4(i)
 \end{eqnarray*}
 with
\begin{eqnarray*}
  I_1(i) & = & \sum_{k=1}^r \frac{[\lambda_k n]}{\hat{n}} \left|  A_{[\lambda_k
             n]}\Phi_i(T^{a_k}y) -  A_{[\lambda_k n]} \Phi_i(x_k)
      \right|
\\
I_2(i) & = &  \sum_{k=1}^r \frac{[\lambda_k n]}{\hat{n}} \left|
A_{[\lambda_k n]} \Phi_i(x_k) -   \int \Phi_i d \nu_k
      \right|\\
I_3(i) & = &  \sum_{k=1}^r  \left| \frac{[\lambda_k n]}{\hat{n}} -
 \lambda_k \right| \int |\Phi_i| d \nu_k
      \\
I_4(i) & = & \frac{1}{\hat{n}} \sum_{k=2}^r \sum_{j=a_k-m}^{a_k
-1} |\Phi_i(T^j
         y)|.
\end{eqnarray*}
Since $[\lambda_k n]\le \lambda_k\hat{n}$ and $x_k $ satisfies
(\ref{7}), by (\ref{3}) we get
\begin{eqnarray*}
  \sum_{i=1}^\infty p_i I_1(i) & \le & \sum_{i=1}^\infty p_i \sum_{k=1}^r \lambda_k  \mbox{\rm var} (\Phi_i, \epsilon)
  =\sum_{i=1}^\infty p_i   \mbox{\rm var} (\Phi_i, \epsilon)
  <\delta.
\end{eqnarray*}
Since $x_k \in Y_k(N)$ and $[\lambda_k n]\ge N$, we have
\begin{eqnarray*}
 \sum_{i=1}^\infty p_i I_2(i) & \le &  \sum_{i=1}^\infty p_i \sum_{k=1}^r \lambda_k  \left|
A_{[\lambda_k n]} \Phi_i(x_k) -   \int \Phi_i d \nu_k
      \right| \\
      & =&
      \sum_{k=1}^r \lambda_k  \sum_{i=1}^\infty p_i  \left|
A_{[\lambda_k n]} \Phi_i(x_k) -   \int \Phi_i d \nu_k
      \right| \le \delta \sum_{k=1}^r  \lambda_k  = \delta.
\end{eqnarray*}
Since $\|\Phi_i\|\le 1$, we have
      \begin{eqnarray*}
\sum_{i=1}^\infty p_i  I_3(i) & \le &  \sum_{k=1}^r \left|
\frac{[\lambda_k n]}{\hat{n}} -
 \lambda_k \right|<\delta
      \\
\sum_{i=1}^\infty p_i I_4(i) & \le & \frac{(r-1)
m(\epsilon)}{\hat{n}} \sum_{i=1}^\infty p_i
 = \frac{(r-1) m(\epsilon)}{\hat{n}}<\delta
\end{eqnarray*}
when $n$ is sufficiently large because $\hat{n} \to \infty$ and
$\frac{[\lambda_k n]}{\hat{n}}\to \lambda_k$.

By combining all these estimates, we obtain (\ref{last}). $\Box$

\subsection{Saturatedness of systems with specification}\

In this subsection we will finish our proof of
Theorem~\ref{specification-saturated}, which says that
 systems satisfying the specification property are saturated.
It remains to prove $h_{\rm top}(G_{\mu}) \geq \Lambda_\Phi(\a)$.
In fact, by Proposition \ref{lamb>hu}, we will have $h_{\rm
top}(G_\mu) \ge h_\mu$. On the other hand, it was known to Bowen
\cite{Bowen73} that $h_\mu \ge h_{\rm top}(G_\mu)$. So, we will
get $h_{\rm top}(G_\mu) = h_\mu$.

\begin{proposition}\label{htop>lamb}
$h_{\rm top}(G_{\mu}) \geq \Lambda_\Phi(\a)$
\end{proposition}

\begin{proof}
It suffices to prove $h_{\rm top}(G_{\mu}) \geq \Lambda_\Phi(\a)
-\r$ for any $\r>0$. To this end, we will construct dynamical
Moran subsets of $G_\mu = X_\Phi(\alpha)$, which approach
$X_\Phi(\alpha)$. The construction is based on separated sets of
$X_\Phi (\a,\d,n)$.

Let $m_k= m(2^{-k}\e)$  be the constants in the definition of
specification. By the definition of $\Lambda_\Phi(\alpha)$ (see
(\ref{Lambda})), when $\epsilon>0 $ is small enough there exist a
sequence of positive numbers $\{\d_k\}$ decreasing to zero and a
sequence of integers $\{n_k\}$ increasing to the infinity such
that
\begin{equation}\label{nk}
     n_k \geq 2^{m_k}
\end{equation}
and that for any $k\ge 1$ we can find a $(n_k, 5\e)$-separated set
$W_k$ of $X_\Phi(\alpha, \d_k, n_k)$ with
\begin{equation}\label{Mk}
      M_k:= \sharp W_k \ge \exp(n_k(\Lambda_\Phi (\a) - \r)).
\end{equation}
Choose a sequence of integers  $\{N_k\}$ such that
\begin{equation}\label{Nk-1}
   N_1=1
\end{equation}
\begin{equation}\label{Nk-2}
   N_k \geq 2^{n_{k+1} + m_{k+1}}, \qquad k\geq 2
\end{equation}
\begin{equation}\label{Nk-3}
   N_{k+1} \geq 2^{N_1n_1+N_2(n_2+m_2) \cdots +N_k(n_k+ m_k)},  \qquad k \geq 1
\end{equation}
Consider  the dynamical Moran set
$F=F(\e,\{W_k\},\{n_k\},\{N_k\})$ as we constructed in the last
subsection. From (\ref{nk})-(\ref{Nk-3}), we get
$$
   \liminf_{n \to \infty} \frac{1}{n}(N_1 \log M_1
      +\cdots +N_k \log M_k + p \log M_{k+1}) \geq
      \Lambda_\Phi(\a) - \r.
$$
By Proposition \ref{entofF}, we have
$$
    h_{\rm top}(F) \geq \Lambda_\Phi(\a)- \r.
$$
Thus we have only to prove $
    F \subseteq X_\Phi (\a).
$ Or equivalently
\begin{equation}\label{FsubsetX}
    \lim_{n\to \infty} \frac{1}{n} \sum_{i=0}^{\infty} p_i \left| S_n \Phi_i(x)
    - n
    \a_i \right| =0   \qquad (x \in F).
\end{equation}

Let us use the same notations as in the last subsection including
$\ell_k$ , $t_k$ , $D_k$ and $L_k$ etc.

Fix $n\ge 1$.  Let $k\ge 1 $ and $0\le p < N_{k+1}$ be the
integers, which depend on $n$,  such that
$$
    \ell_k + p(m_{k+1}+ n_{k+1}) < n \leq \ell_k + (p+1)(m_{k+1}+ n_{k+1})
$$
Write
$$
   q=n-\big(\ell_k +p(m_{k+1} + n_{k+1}) \big),   \qquad
   b_s=(s-1)(m_{k+1}+n_{k+1}).
$$
Decompose the interval $[0,n)$ ($\subset\mathbb{N}$) into small
intervals
 $$
 [0,n) = [0, \ell_k) \bigcup [\ell_k, \ell_k+p(m_{k+1}+n_{k+1}))
 \bigcup [\ell_k+p(m_{k+1}+n_{k+1}), n)
  $$
and decompose still $[\ell_k, \ell_k+p(m_{k+1}+n_{k+1}))$ into
intervals alternatively of lengths $n_{k+1}$ and $m_{k+1}$. Then
cut the sum $\sum_{0\le j <n} \Phi_i(T^j x)$ into sums taken over
small intervals.
%
Thus we get
$$
    \left| S_n \Phi_i(x) - n\a_i \right| \leq J_1(i) +J_2(i) +J_3(i) +J_4(i)
$$
where
\begin{eqnarray*}
J_1(i)&=& \left| S_{\ell_k}\Phi_i(x) - \ell_k \a_i\right|\\
J_2(i)&=& \sum_{s=1}^p \left| S_{m_{k+1}}\Phi_i(T^{\ell_k+b_s}x) - m_{k+1} \a_i \right|\\
J_3(i)&=& \sum_{s=1}^p \left| S_{n_{k+1}}\Phi_i(T^{\ell_k+b_s+m_{k+1}}x) - n_{k+1}\a_i \right|\\
J_4(i)&=& \left| S_q\Phi_i(T^{\ell_k+p(m_{k+1}+n_{k+1})}x)- q \a_i
               \right|
\end{eqnarray*}

Since $\|\Phi_i\|\le 1$ (hence $|\alpha_i|\le 1$), we have
$$
    J_2(i) \leq 2\medskip p\medskip m_{k+1},\qquad
    J_4(i) \leq 2\medskip q \leq 2(m_{k+1}+ n_{k+1}).
$$
By (\ref{Nk-2}),  we have
\begin{equation}\label{J2J4}
     \lim_{n \to \infty}\frac{1}{n} \sum_{i=0}^\infty p_i J_2(i)  =
     0,
     \quad
 \lim_{n \to \infty}
   \frac{1}{n} \sum_{i=0}^\infty p_i J_4(i) =0.
\end{equation}

 Now let us deal with $J_1(i)$ and $J_3(i)$.
   We claim that for any $x\in F$
there exists an $\bar{x}\in L_k$ such that
\begin{equation}\label{F-1}
   d_{\ell_k}(\bar{x},x)< \frac{\e}{2^{k-1}},
\end{equation}
and that for all $1\le s \le p$, there exists a point $x_{s}\in
W_{k+1}$ such that
\begin{equation}\label{F-2}
   d_{n_{k+1}}(x_{s}, T^{u_s}x) <
   \frac{\e}{2^{k-1}}
\end{equation}
where $$ u_s = \ell_k + b_s + m_{k+1}.
$$
 In fact, by the
construction of $F$,  there exists a point $z\in L_{k+1}$ such
that
\begin{equation}\label{F-1-1}
  d_{\ell_{k+1}}(z,x) \le \frac{\e}{2^{k}}.
\end{equation}
Assume that $z$ descends from some $\bar{x}\in L_k$ through $y\in
D_{k+1}$. Then
\begin{equation}\label{F-1-2}
   d_{\ell_k}(\bar{x},z)< \frac{\e}{2^{k+1}}
\end{equation}
and
\begin{equation}\label{F-2-1}
   d_{t_{k+1}}(y, T^{\ell_k+ m_{k+1}}z) < \frac{\e}{2^{k+1}}.
\end{equation}
On the other hand,  according to the definition of $D_{k+1}$,
there exists an $x_{s}\in W_{k+1}$ such that
\begin{equation}\label{F-2-2}
  d_{n_{k+1}}(x_{s},T^{b_s}y) < \frac{\e}{2^{k+1}}.
\end{equation}
Now  by the trigonometric inequality, the fact $d_{\ell_k}(z,x)\le
d_{\ell_{k+1}}(z,x)$ and (\ref{F-1-1}) and (\ref{F-1-2}) we get
\begin{eqnarray*}
   d_{\ell_k}(\bar{x},x) \le
                         \frac{\e}{2^{k+1}} + \frac{\e}{2^{k}}
                         <    \frac{\e}{2^{k-1}}.
\end{eqnarray*}
Thus (\ref{F-1}) is proved.  By (\ref{F-1-1}),(\ref{F-2-1}) and
(\ref{F-2-2}), we can similarly prove (\ref{F-1-2}):
\begin{eqnarray*}
    d_{n_{k+1}}(x_{s}, T^{u_s} x)
   &\le&  d_{n_{k+1}}(x_{s}, T^{b_s}y)
          + d_{n_{k+1}}(T^{b_s}y, T^{u_s} z)
          + d_{n_{k+1}}(T^{u_s}z,T^{u_s}x)\\
   &\le&  d_{n_{k+1}}(x_{s},T^{b_s}y) + d_{t_{k+1}}(y,T^{\ell_k+ m_{k+1}}z) + d_{\ell_{k+1}}(z,x) \\
   &<&    \frac{\e}{2^{k+1}} + \frac{\e}{2^{k+1}}  + \frac{\e}{2^{k}}\\
   &=&    \frac{\e}{2^{k-1}}.
\end{eqnarray*}

It is now easy to deal with  $J_3(i)$, which is obviously bounded
by
\begin{eqnarray*} \label{J3-1}
    J_3(i)\le  \sum_{s=1}^p
               \left| S_{n_{k+1}} \Phi_i(T^{u_s} x)-S_{n_{k+1}} \Phi_i(x_{s})
               \right|
         + \sum_{s=1}^p
                \left| S_{n_{k+1}} \Phi_i(x_{s}) -n_{k+1} \a_i
                \right|.
\end{eqnarray*}
Using (\ref{F-2}), we obtain
\begin{eqnarray*} \label{J3-2}
   \left| S_{n_{k+1}} \Phi_i(x_{s}) - S_{n_{k+1}}\Phi_i(T^{u_s}x)
   \right| \leq n_{k+1}{\rm var}(\Phi_i, \e 2^{-(k-1)}).
\end{eqnarray*}
On the other hand, since $
  x_{s} \in W_{k+1} \subseteq X_\Phi(\a,\d_{k+1},n_{k+1}),
$ we have, by definition,
\begin{eqnarray*} \label{J3-3}
    \sum_{i=1}^\infty p_i \left| S_{n_{k+1}}\Phi_i(x_{s}) -
                                 n_{k+1} \a_i \right|
    \leq n_{k+1} \d_{k+1}.
\end{eqnarray*}
Then, combining the last three estimates, and using the facts
$\sum_{j=1}^\infty p_j =1$ and $p n_{k+1}\le n$, we get
$$
   \frac{1}{n}\sum_{i=1}^\infty p_i J_3(i) \leq
   \sum_{i=1}^\infty p_i{\rm var}(\Phi_i,
   \e 2^{-(k-1)}) + \d_{k+1}.
$$
Since $k$ can be arbitrarily large, we finally get
\begin{equation}\label{J3}
    \lim_{n\to \infty} \frac{1}{n}\sum_{i=1}^\infty p_i J_3(i) =0.
\end{equation}

Now it remains to prove
\begin{eqnarray}\label{J1}
    \lim_{n \to \infty}  \frac{1}{n}\sum_{i=1}^\infty p_i J_1(i)=0.
\end{eqnarray}
Observe that
$$
   J_1(i)\leq \left| S_{\ell_k}\Phi_i(x) - S_{\ell_k}\Phi_i(\bar{x})\right| +
                 \left| S_{\ell_k}\Phi_i(\bar{x}) - \ell_k \a_i
                 \right|.
$$
 By (\ref{F-1}), we have $|S_{\ell_k}\Phi_i(x) - S_{\ell_k}\Phi_i(\bar{x})|\le \ell_k {\rm var}(\Phi_i, \e
 2^{-(k-1)})$. Then
\begin{eqnarray*}
    J_1(i)
         \leq  \ell_k {\rm var}(\Phi_i, \e 2^{-(k-1)}) + R_{k,i}
\end{eqnarray*}
where
$$
     R_{k,i}= \max_{z\in L_k} \left| S_{\ell_k}\Phi_i(z) -
    \ell_k \a_i \right|.
$$
Since ${\rm var}(\Phi_i, \e 2^{-(k-1)})$ tends to zero as $k \to
\infty$, the desired claim (\ref{J1}) is reduced to
\begin{eqnarray}\label{R}
    \lim_{n \to \infty} \frac{1}{n} \sum_{i=1}^\infty p_i R_{k,i}=0.
\end{eqnarray}
We need two lemmas to estimate $R_{k, i}$.

\begin{lemma}\label{Dk}
For any $y \in D_{k+1}$, we have
\begin{eqnarray*}
   & & \sum_{i=1}^\infty p_i \left| S_{t_{k+1}}\Phi_i(y) -t_{k+1} \a_i \right|\\
   &\leq&
    \sum_{i=1}^\infty p_i
    N_{k+1}n_{k+1}{\rm var}(\Phi_i,\e2^{-(k+1)}) + 2(N_{k+1}-1)m_{k+1} +
     N_{k+1} n_{k+1} \d_{k+1}.
\end{eqnarray*}
\end{lemma}

\begin{proof}
For  any $s=1,\ldots,N_{k+1}$, there exists $x_{s}\in W_{k+1}$
such that
\begin{equation}\label{dnk+1}
   d_{n_{k+1}}(x_{s},T^{b_s}y) < \frac{\e}{2^{k+1}}
\end{equation}
where $
   b_s= (s-1)(m_{k+1}+n_{k+1}).
$ Write
$$
   S_{t_{k+1}}\Phi_i(y)= \sum_{s=1}^{N_{k+1}}
   S_{n_{k+1}}\Phi_i(T^{b_s}y)
   + \sum_{s=1}^{N_{k+1}-1}S_{m_{k+1}}\Phi_i(T^{b_s+n_{k+1}}y).
$$
Then
\begin{eqnarray*}
   & &    \left| S_{t_{k+1}}\Phi_i(y)- t_{k+1}\a_i \right|\\
   &\le&  \sum_{s=1}^{N_{k+1}} \left|
            S_{n_{k+1}}\Phi_i(T^{b_s}y) - n_{k+1}\a_i \right|
   + \sum_{s=1}^{N_{k+1}-1} \left|
            S_{m_{k+1}}\Phi_i(T^{b_s+n_{k+1}}y)- m_{k+1}\a_i
            \right|.
\end{eqnarray*}
Since $
   x_{s}\in W_{k+1}\subseteq X_\Phi(\a,\d_{k+1},n_{k+1}),
$ by (\ref{dnk+1}), we have
\begin{eqnarray*}
    & &    \sum_{i=1}^{\infty}p_i \left| S_{n_{k+1}}\Phi_i(T^{b_s}y) - n_{k+1}\a_i \right|\\
    &\leq& \sum_{i=1}^{\infty}p_i\left| S_{n_{k+1}}\Phi_i(T^{b_s}y) - S_{n_{k+1}}\Phi_i(x_s) \right|
          + \sum_{i=1}^{\infty}p_i\left| S_{n_{k+1}}\Phi_i(x_s) - n_{k+1}\a_i \right|\\
    &\leq& \sum_{i=1}^{\infty}p_i n_{k+1}{\rm var}(\Phi_i,\e2^{-(k+1)})+ n_{k+1}
    \d_{k+1}.
\end{eqnarray*}
On the other hand,
$$
  \left| S_{m_{k+1}}\Phi_i(T^{b_s+n_{k+1}}y)- m_{k+1}\a_i \right|
  \le 2m_{k+1}.
$$
Now it is easy to conclude.
\end{proof}

\begin{lemma}\label{Rki}
$$
   \sum_{i=1}^\infty p_i R_{k,i}\leq 2\sum_{i=1}^\infty p_i \sum_{j=1}^k
   \ell_j {\rm var}(\Phi_i,\e 2^{-j})+ 2\sum_{j=1}^k N_j m_j + \sum_{j=1}^k \ell_j
   \d_j.
$$
\end{lemma}

\begin{proof}
We prove it by induction on $k$. When $k=1$, we have $
   L_1=D_1=W_1 \subseteq X_\Phi(\a,\d_1,n_1)
$ and then
$$
     \sum_{i=1}^\infty p_i  R_{1,i}\leq n_1\d_1= \ell_1\d_1.
$$
Suppose the lemma holds for $k$. For any $z\in L_{k+1}$ there
exist $x\in L_k$ and $y\in D_{k+1}$, such that
$$
    d_{\ell_k}(x,z) < \frac{\e}{2^{k+1}}, \qquad d_{t_{k+1}}(y, T^{\ell_k +
    m_{k+1}}z) < \frac{\e}{2^{k+1}}.
$$
Write
$$
   S_{\ell_{k+1}}\Phi_i(z) = S_{\ell_{k}}\Phi_i(z)+ S_{m_{k+1}}\Phi_i(T^{\ell_k}z) +
   S_{t_{k+1}}\Phi_i(T^{\ell_k+ m_{k+1}}z).
$$
Then $ \left| S_{\ell_{k+1}} \Phi_i (z) -\ell_{k+1} \a_i\right|$
is bounded by
\begin{eqnarray*}
      \left| S_{\ell_k} \Phi_i (z) - \ell_k \a_i\right|
          +  \left| S_{m_{k+1}} \Phi_i(T^{\ell_k} z) - m_{k+1}\a_i\right|
           +  \left| S_{t_{k+1}} \Phi_i (T^{\ell_k+m_{k+1}} z) -t_{k+1}
           \a_i\right|.
\end{eqnarray*}
Notice that
\begin{eqnarray*}
   \left| S_{\ell_k} \Phi_i (z) - \ell_k \a_i\right|
   &\leq& \left| S_{\ell_k} \Phi_i (z) - S_{\ell_k} \Phi_i (x) \right|
             +  \left| S_{\ell_k} \Phi_i (x) - \ell_k \a_i\right|\\
   &\le&  \ell_k {\rm var}(\Phi_i,\e 2^{-(k+1}) + R_{k,i},
\end{eqnarray*}
\begin{eqnarray*}
      \left| S_{m_{k+1}} \Phi_i(T^{\ell_k} z) - m_{k+1}\a_i\right|
      & \le&  2m_{k+1}
\end{eqnarray*}
and
\begin{eqnarray*}
    & &    \left| S_{t_{k+1}} \Phi_i (T^{\ell_k+m_{k+1}} z) -t_{k+1} \a_i\right|\\
    &\le&  \left| S_{t_{k+1}} \Phi_i (T^{\ell_k+m_{k+1}} z) -S_{t_{k+1}} \Phi_i(y) \right|
        +  \left| S_{t_{k+1}} \Phi_i (y) -t_{k+1} \a_i\right|\\
    &\le&  t_{k+1} {\rm var}(\Phi_i,\e 2^{-(k+1)}) +
           \left| S_{t_{k+1}} \Phi_i (y) -t_{k+1} \a_i\right|.
\end{eqnarray*}
 By Lemma~\ref{Dk}, we have
\begin{eqnarray*}
   \sum_{i=1}^\infty p_i R_{k+1,i} &\le& \sum_{i=1}^\infty p_i R_{k,i}
                    +\sum_{i=1}^\infty p_i (\ell_k + t_{k+1} + N_{k+1}n_{k+1})
                     {\rm var}(\Phi_i,\e/2^{k+1}) \\
       & &    + 2N_{k+1}m_{k+1} +
                 N_{k+1}n_{k+1}\d_{k+1}.
\end{eqnarray*}
Then according to the induction hypothesis the Lemma holds for
$k+1$, because
$$
    \ell_k+ t_{k+1} \leq \ell_{k+1}, \quad  N_{k+1}n_{k+1} \le
    \ell_{k+1}.
$$
\end{proof}

Let us finish the proof of Theorem~\ref{specification-saturated}
by showing (\ref{R}). Since $n_j \geq 2^{m_j}$, we have
$$
    \frac{N_jm_j}{\ell_j} \leq \frac{N_jm_j}{N_j(n_j+m_j)}=
    \frac{m_j}{n_j+ m_j} \rightarrow 0 \quad (j\to \infty).
$$
Then the estimate in Lemma~\ref{Rki} can be written as
$$
    \sum_{i=1}^\infty p_i R_{k,i} \leq \sum_{j=1}^k \ell_j c_j
$$
where $c_j \to 0$ ($j\to \infty$). By (\ref{Nk-3}), we have $
    \ell_k\geq 2^{\ell_{k-1}}.
$ It follows that
$$
    \frac{1}{\ell_k}
       \sum_{i=1}^{\infty } p_i R_{k,i}  \leq c_k +
    \frac{1}{\ell_k} \sum_{i=1}^{k-1} c_j \ell_j.
$$
This implies (\ref{R}).
\end{proof}

\setcounter{equation}{0}
\section{Variational principle}

In this section, we prove variational principles for saturated
systems (Theorem~\ref{VP} and Theorem~\ref{saturated}).

\subsection{Proofs of Theorems \ref{VP} and \ref{saturated}}\

The proof of Theorem~\ref{saturated} is similar to that of
Theorem~\ref{VP}.

\textit{Proof of Theorem~\ref{VP} (a).\/} It suffices to prove
that if there exists a point $x \in X_\Phi(\alpha; E)$, then $
\mathcal{M}_\Phi(\alpha; E)\not= \emptyset$. That $x \in
X_\Phi(\alpha; E)$ means
\begin{equation}\label{3A1}
      \limsup_{n \to \infty} \langle A_n \Phi(x), w\rangle \le \langle
      \alpha, w \rangle \qquad (\forall w \in E).
\end{equation}
Let $\mu$ be a weak limit of $n^{-1}\sum_{j=0}^{n-1}
\delta_{T^jx}$. That is to say, there exists a sequence $n_m$ such
that
\begin{equation}\label{3A2}
        \lim_{m \to \infty} \frac{1}{n_m}
             \sum_{j=0}^{n_m-1} f(T^j x) = \int f d\mu
\end{equation}
for all scalar continuous functions $f$. We deduce from
(\ref{3A1}) and (\ref{3A2}) that for all $w \in E$ we have
$$
     \int \langle \Phi, w \rangle d\mu
     = \lim_{m \to \infty} \langle A_{n_m} \Phi(x), w\rangle
             \le
             \limsup_{n \to \infty} \langle A_n \Phi(x), w\rangle \le \langle \alpha, w
      \rangle.
$$
So $\mu \in \mathcal{M}_\Phi(\alpha; E)$.

\textit{Proof of Theorem~\ref{VP} (b).\/} Let $t = \sup_{\mu \in
\mathcal{M}_\Phi(\alpha; E)} h_\mu$.  What we have just proved
above may be stated as follows: if $x \in X_\Phi(\alpha; E)$, then
$$V(x) \subset \mathcal{M}_\Phi(\alpha; E).$$
It follows that    $h_\mu \le t$ for any $\mu \in V(x)$. Thus
\begin{eqnarray*}
     X_\Phi(\alpha; E)
     & \subset &
        \left \{ x \in X: \ \forall \ \mu \in V(x)\
       \mbox{\rm satisfying } \ h_\mu \le t\right\}\\
       & \subset & \left \{ x \in X: \ \exists \ \mu \in V(x)\
       \mbox{\rm satisfying } \ h_\mu \le t\right\}.
\end{eqnarray*}
Then, due to Lemma~\ref{Bowen}, we get $h_{\rm top}
(X_\Phi(\alpha; E)) \le t$.

Now we prove the converse inequality. For any $\mu
\in\mathcal{M}_\Phi(\alpha; E)$, consider $G_\mu$ the set of
generic points. We have
$$
      G_\mu \subset X_\Phi(\alpha; \mathbb{B}) \subset X_\Phi(\alpha; E).
$$
The second inclusion is obvious and the first one is a consequence
of the fact that $x \in G_\mu$ implies $\lim_{n \to \infty}
A_n\Phi(x) = \int \Phi d \mu =\alpha$ in the $\sigma(\mathbb{B}^*,
\mathbb{B})$-topology. Thus
$$
    h_{\rm top} (X_\Phi(\alpha; E)) \ge h_{\rm top} (G_\mu).
$$
Since $\mu $ is an arbitrary invariant measure in
$\mathcal{M}_\Phi(\alpha; E)$, we can finish the proof because the
system $(X, T)$ is saturated (i.e. $h_{\rm top} (G_\mu) =h_\mu$).
$\Box$
\medskip

It is useful to point out the following facts appearing in the proof:\\
\indent (i)\ \ If $x \in X_\Phi(\alpha, E)$, then
    $V(x) \subset \mathcal{M}_\Phi(\alpha, E)$.\\
\indent (ii) \ We have
$$
      \bigcup_{\mu \in \mathcal{M}_\Phi(\alpha, E) }
      G_\mu
      \subset X_\Phi(\alpha, E)
      \subset \bigcup_{\mu \in \mathcal{M}_\Phi(\alpha, E) }
      \widetilde{G}_\mu
      $$
      with $\widetilde{G}_\mu = \{x\in X: V(x) \ni \mu \}$.  It is
worth to notice the fact that all the $G_\mu$ are disjoint.

It is clear that $\mathcal{M}_\Phi(\alpha, E)$ is a compact convex
subset of the space $\mathcal{M}_{\rm inv}$ of Borel probability
invariant measures. If $h_\mu$, as a function of $\mu$, is upper
semi-continuous on $\mathcal{M}_{\rm inv}$, then the supremum in
the variational principle is attained by some invariant measure,
called the maximal entropy measure in $\mathcal{M}_\Phi(\alpha,
E)$. Usually, the structure of $\mathcal{M}_\Phi(\alpha, E)$ is
complicated. But it is sometimes possible to calculate the maximal
entropy.
\medskip

\textit{Proof of Theorem~\ref{saturated}.\/} Let $x \in
X_\Phi^\Psi(\beta)$ and let $\mu$ be a weak limit of
$n^{-1}\sum_{j=0}^{n-1} \delta_{T^j x}$. Then there exists a
subsequence of integers $\{n_m\}$ such that $A_{n_m}\Phi(x)$ tends
to $\int \Phi d \mu$ in the weak star topology as $m \to \infty$
because we have an expression similar to (\ref{3A2}) with $f =
\langle \Phi, w \rangle$ ($w \in \mathbb{B}$ being arbitrary).
Hence
$$
    \Psi \left( \int \Phi d \mu \right)
     = \lim_{ m\to \infty} \Psi(A_{n_m}\Phi(x))
     = \lim_{ n\to \infty} \Psi(A_{n}\Phi(x)) = \beta.
$$
Thus we have proved that $\mu \in \mathcal{M}_\Phi^\Psi(\beta)$.
That is to say
$$
         V(x) \subset  \mathcal{M}_\Phi^\Psi(\beta) \quad (\forall x \in
         X_\Phi^\Psi(\beta)).
$$
It follows that (a) holds and that due to Lemma~\ref{Bowen} we
have
$$
      h_{\rm top} (X_\Phi^\Psi(\beta)) \le \sup_{\mu \in
      \mathcal{M}_\Phi^\Psi(\beta)} h_\mu.
$$

The converse inequality is a consequence of the variational
principle (\ref{variational-principle}) and the relationship
$$
    X_\Phi^\Psi(\beta) \supset \widehat{X}_\Phi^\Psi(\beta)
    =\bigcup_{\alpha: \Psi(\alpha)=\beta} \ \  X_{\Phi}(\alpha).
$$
In fact,
\begin{eqnarray*}
      h_{\rm top}   (X_\Phi^\Psi(\beta))\ge  h_{\rm top}
      (\widehat{X}_\Phi^\Psi(\beta))
         & \ge & \sup_{\alpha: \Psi(\alpha)=\beta} \ \  h_{\rm top}
         (X_{\Phi}(\alpha))\\
         &= & \sup_{\alpha: \Psi(\alpha)=\beta} \ \  \sup_{\mu \in \mathcal{M}_\Phi(\alpha)}
         h_\mu\\
         & = & \sup_{\mu \in \mathcal{M}_\Phi^\Psi(\beta)}
         h_\mu.
\end{eqnarray*}

\subsection{$\ell^\infty(\mathbb{Z})$-valued ergodic average}\

Let us consider the special case where $\mathbb{B}=
\ell^1(\mathbb{Z})$.  Then $\mathbb{B}^* =
\ell^\infty(\mathbb{Z})$. Any $\ell^\infty(\mathbb{Z})$-valued
function $\Phi$ can be written as
$$
  \Phi(x) = \bigl(\Phi_n(x)\bigr)_{n\in{\mathbb
 Z}},\quad\text{with}\quad \sup_n |\Phi_n(x)| <\infty.
$$
Recall that $\ell^\infty(\mathbb{Z})$ is equipped with the
$\sigma(\ell^\infty, \ell^1)$-topology.  A
$\ell^\infty(\mathbb{Z})$-valued function $\Phi$ is continuous if
and only if all coordinate functions $\Phi_n : X \to \mathbb{R}$
are continuous, because for any $w=(w_n)_{n\in Z} \in \ell^1$ we
have
$$
   \langle \Phi(x), w \rangle = \sum_{n\in Z} w_n \Phi_n(x).
$$

Let us give an application of the variational principle in this
setting. Let $I$ be a finite or infinite subset of positive
integers. Let $\{\Phi_i\}_{i\in I}$ be a family of real continuous
functions defined on $X$.  We suppose that $\sup_{i\in I}
\|\Phi_i\|_{C(X)}<\infty$.  For two given sequences of real
numbers $ \mathbf{a}=\{a_i\}_{i\in I}$ and
$\mathbf{b}=\{b_i\}_{i\in I}$, we denote by $S(\mathbf{a},
\mathbf{b})$ the set of points $x \in X$ such that
$$
         a_i \le \liminf_{n\to \infty} A_n\Phi_i(x)
         \le \limsup_{n\to \infty} A_n\Phi_i(x) \le b_i
         \qquad (\forall i \in I).
$$

\begin{corollary} Suppose that the system $(X, T)$ is saturated.
The topological entropy of $S(\mathbf{a}, \mathbf{b})$ defined
above is equal to the supremum of the measure-theoretical
entropies $h_\mu$ for those invariant measures $\mu$ such that
$$
         a_i \le \int \Phi_i d \mu  \le b_i
         \qquad (\forall i \in I).
$$
\end{corollary}

\proof For $n \in \mathbb{Z}$, let $ e_n $ be the $n^\text{th}$
element of the canonical basis of $\ell^1(\mathbb{Z})$. Let $\Phi$
be a function whose $n$th coordinate and $-n^\text{th}$ coordinate
are respectively equal to $\Phi_{n}$ and $-\Phi_n$ for each $n \in
I$ (other coordinates may be taken to be zero). Take the set $E
\subset \ell^1$, which consists of $e_i$ and $e_{-i}$ for $i \in
I$. Take $\alpha \in \ell^\infty$ such that $\alpha_{-i} = - a_i$
and $\alpha_i = b_i$ for $i \in I$. Now we can directly apply the
variational principle by noticing that
$$
     \langle \Phi,e_i \rangle =  \Phi_i,
     \qquad  \langle \Phi,e_{-i} \rangle
     = -\Phi_i  \qquad \qquad (i \in I).
$$
$\Box$

The result contained in this corollary is new, even when $I$ is
finite.  If $I$ is finite and if $a_i=b_i$ for $i\in I$, the
preceding corollary allows one to recover the results in
\cite{FFW} and \cite{TV}.

The validity of the variational principle is to some extent
equivalent to the fact that the system $(X, T)$ is saturated.

\begin{theorem} Let $(X, T)$ be a compact dynamical system.
The system is saturated if and only if the variational principle
(Theorem~{\ref{VP}} (b)) holds for all real Banach spaces
$\mathbb{B}$.
\end{theorem}

\proof It remains to prove that the variational principle implies
the saturation of the system.

Take a countable  set $\{f_i\}_{i \in \mathbb{N}}$ which is dense
in the unit ball of $C(X)$ ($C(X)$ being the space of all real
valued continuous functions on $X$). Consider the function
$$\Phi=(f_1, f_2, \cdots, f_n, \cdots)
$$
which  takes values in $\mathbb{B}=\ell^\infty(\mathbb{N})$. For
any invariant measure $\mu \in \mathcal{M}_{\rm inv}$, define
$$
    \alpha = \left(\int f_1 d\mu, \int f_2 d\mu, \cdots\right) \in \ell^\infty(\mathbb{N}).
$$
It is clear that $\mathcal{M}_\Phi(\alpha) = \{\mu\}$. Then the
variational principle implies $h_{\rm top} (X_\Phi(\alpha)) =
h_\mu$. This finishes the proof because $X_\Phi(\alpha)$ is
nothing but $G_\mu$. $\Box$

\section{An example: recurrence in an infinite number of cylinders}

We have  got a general variational principle. In order to apply
this principle, one of the main questions is to get information
about the convex set $\mathcal{M}_\Phi(\alpha, E)$ and the maximal
entropy measures contained in it and to compute the maximal
entropy. Let us consider the symbolic dynamical system $(\{0,
1\}^\mathbb{N}, T)$, $T$ being the shift. The structure of the
space $\mathcal{M}_{\rm inv}$ is relatively simple.  To illustrate
the main result, we shall consider a special problem of recurrence
into a countable number of cylinders.

\subsection{Symbolic space}\

Let $X=\{0, 1\}^\mathbb{N}$ and $T$ be the shift transformation.
As
 usual, an $n$-cylinder in $X$ determined by a word $w=x_1x_2\cdots
 x_n$ is denoted by $[w]$ or $[x_1, \cdots, x_n]$. For any word $w$,
 define the recurrence to the cylinder $[w]$ of $x$ by
 $$
        R(x, [w]) = \lim_{n\to \infty} \frac{1}{n}\sum_{j=0}^{n-1}
            1_{[w]}(T^j x)
 $$
 if the limit exists.

Let $\mathcal{W}=\{w_i\}_{i \in I}$ with $I \subset \mathbb{N}$ be
a finite or infinite set of words.  Let $\alpha =\{a_i\}_{i \in
I}$ be a sequence of non-negative numbers.  We are interested in
the following recurrence set
$$
    R(\mathbf{a};\mathcal{W}) = \left\{x \in X: R(x, [w_i]) = a_i
    \ \mbox{\rm for \ all } \ i \in I\right\}
$$
whose topological entropy will be computed by the variational
principle which takes a simpler form.

\begin{corollary}\label{Cor-Example} We have $h_{\rm top}(R(\mathbf{a};\mathcal{W}))=
\max_{\mu \in \mathcal{M} (\mathbf{a, \mathcal{W}})} h_\mu$ where
$$
     \mathcal{M} (\mathbf{a, \mathcal{W}})
       = \left\{\mu \in \mathcal{M}_{\rm inv}:   \mu([w_i]) = a_i
        \ \mbox{\rm for \ all } \ i \in I\right\}.
$$
\end{corollary}
Remark that the shift transformation on the symbolic space is
expansive. Thus the entropy function $h_{\mu} $ is upper
semi-continuous(\cite{Walters}, p. 184). Hence we can obtain the
supremum in Corollary \ref{Cor-Example}.

Recall that
     any Borel probability measure on $X$ is uniquely determined by its values on
     cylinders. In other words, a function $\mu$ defined on all
     cylinders can be extended to be a Borel probability measure if and only if
         $$
         \sum_{x_1,\cdots,x_n}\mu([x_1,\cdots,x_n])=1
         $$
     and
     $$
        \sum_{\epsilon \in \{0, 1\}} \mu([x_1,\cdots,x_{n-1}, \epsilon])
        = \mu([x_1,\cdots,x_{n-1}]).
     $$
      Such a probability measure $\mu$ is invariant if and only if
    $$
    \sum_{\epsilon\in \{0, 1\} }\mu([\epsilon,x_2,\cdots,x_n])=\mu([x_2,\cdots,x_n]).
    $$

    The entropy $h_{\mu}$ of any invariant measure
    $\mu\in\mathcal{M}_{\rm inv}$ can be expressed as follows
   $$h_{\mu}=\lim_{n\to\infty}\sum_{x_1,\cdots,x_n}-\mu([x_1,\cdots,x_n])
   \log \frac{\mu([x_1,\cdots,x_n])}{\mu([x_1,\cdots,x_{n-1}])}.
   $$
 The sum in the above expression which we will denote by $h_{\mu}^{(n)} $ is nothing but a conditional entropy
 of $\mu$ and it is also the entropy of an $(n-1)$-Markov measure
 $\mu_n$, which tends towards $\mu$ as $n$ goes to $\infty$.

 A Markov measure of order $k$ is an invariant measure
  $\nu\in\mathcal{M}_{\rm inv}$ having the following Markov property:
  for all $n>k$ and all $(x_1, \cdots, x_n)\in \{0, 1\}^n$
   $$\frac{\nu([x_1,\cdots,x_n])}{\nu([x_1,\cdots,x_{n-1}])}
=\frac{\nu([x_{n-k},\cdots,x_n])}{\nu([x_{n-k},\cdots,x_{n-1}])}.$$
A Markov measure of order $k$ is uniquely determined by its values
on the $(k+1)$-cylinders.  The preceding approximating Markov
measure $\mu_n$ has the same values as $\mu$ on $n$-cylinders.

To apply the above corollary, we have to maximize the entropy
$h_\mu$ among all invariant measures $\mu$ with constraints
$\mu([w_i]) = a_i$ for $i \in I$. The entropy $h_\mu$ is a
function of an infinite number of variables $\mu([w])$. So we have
to maximize a function of an infinite number of variables.
However, in some cases it suffices to reduce the problem to
maximize the conditional entropy which is a function of a finite
number of variables.

Denote by $|w| $ the length of the word $w$. Let $\mathcal W_n :=
\{ w \in \mathcal W: |w| \le n \} $ and $\mathcal
M(\mathbf{a},\mathcal W _n) := \{ \mu \in \mathcal M_{\rm inv}:
\mu([w_i])=a_i, w_i \in \mathcal W_n \}$. Let  $\mu^* $ be a
maximal entropy measure over $\mathcal M(\mathbf{a},\mathcal W)$
and $\mu_n^* $ be the $(n-1) $-Markov measure which converges to
$\mu^* $. Let  $\mu^{(n)} $ be a maximal entropy measure over
$\mathcal M(\mathbf{a},\mathcal W _n)$. Then
\begin{eqnarray*}
 h_{\mu^*} = \lim_{n\to \infty} h_{\mu_n^*} \le \liminf_{n\to \infty}
 h_{\mu^{(n)}} \le \limsup_{n\to \infty}
 h_{\mu^{(n)}} \le h_{\mu^*}.
\end{eqnarray*}
Hence
\begin{eqnarray*}
\lim_{n\to \infty} h_{\mu^{(n)}} = h_{\mu^*}= \max_{\mu \in
\mathcal{M} (\mathbf{a, \mathcal{W}})} h_\mu.
\end{eqnarray*}
However, for any measure $\mu \in \mathcal M(\mathbf{a},\mathcal W
_n)  $, we have $ h_{\mu}=h_{\mu_n}=h_{\mu}^{(n)}$, where $\mu_n $
is the $(n-1) $-Markov measure which converges to $\mu$
(\cite{FFW}). Thus, $\mu^{(n)} $ is the maximal point of the
conditional entropy function $h_{\mu}^{(n)} $.

Thus we have proved the following proposition.

\begin{proposition}\label{approximation}
The maximal entropy over $\mathcal M(\mathbf{a},\mathcal W)$ can
be approximated by the maximal entropies over $\mathcal
M(\mathbf{a},\mathcal W_n)$'s.
\end{proposition}

\setcounter{equation}{0} \subsection{Example: Frequency of dyadic
digital blocks}\

Let us consider a special example:
$$
      \mathcal{W} = \{[0], [0^2], \cdots, [0^n], \cdots\}
$$
where $0^k$ means the word with $0$ repeated $k$ times.

\begin{theorem} \label{Frequency}
 Let $\mathcal{W} = \{[0^n]\}_{n\ge 1}$ and
$\mathbf{a} = \{a_n\}_{n\ge 1} \subset \mathbb{R}^+$. We have
   \\
   \indent {\rm (a)} \  $R(\mathbf{a}; \{[0^n]\}_{n\ge 1} )\not= \emptyset$
   if and only if
   \begin{equation}
    1=a_0 \ge a_1 \ge a_2 \ge \dots; \quad a_i - 2 a_{i+1} +
    a_{i+2}\ge 0  \ \ (i\ge 0).
\end{equation}
\indent {\rm (b)} If the above condition is satisfied, then
\begin{equation}
   h_{\rm top} (R(\mathbf{a}; \mathcal{W} ))
   = -h(1-a_1) + \sum_{i=0}^\infty h( a_i - 2 a_{i+1} +
    a_{i+2})
\end{equation}
where $h(x) = -x \log x$.
\end{theorem}

The proof of the above theorem is decomposed into several lemmas
which actually allow us to find the unique invariant measure of
maximal entropy and to compute its entropy.

Let $\mu$ be an invariant measure.  The consistence and the
invariance of the measure imply that we may partition all
$(n+2)$-cylinders into groups of the form
      $$ \{[0w0], [0w1], [1w0], [1w1] \} $$
such that  the measures $\mu([0w0]), \mu([0w1]),
\mu([1w0]),\mu([1w1])$ are linked each other through the measures
$\mu([0w]), \mu([w0]), \mu([w1]), \mu([1w])$ of $(n+1)$-cylinders.
More precisely, if we write $p_w = \mu([w])$, then for any word
$w$ of length $n$, we have
\begin{eqnarray*}
     p_{0w0}+p_{0w1}
     &= &p_{0w}\\
      p_{1w0}+p_{1w1}& =& p_{1w}\\
     p_{0w0}+p_{1w0} & =& p_{w0}\\
       p_{0w1}+p_{1w1} &=& p_{w1}
\end{eqnarray*}

\begin{lemma}\label{*0^k*}
Suppose $\mu \in \mathcal{M}(\mathbf{a}, \mathcal{W})$. If $w=0^n$
with $n\ge 0$, we have
\begin{eqnarray}
   p_{00^{n}0} &=& a_{n+2} \label{000}\\
   p_{00^{n}1}& =& a_{n+1}-a_{n+2}\label{001}\\
   p_{10^{n}0}& =& a_{n+1}-a_{n+2} \label{100}\\
    p_{10^n1}&=&
   a_n-2a_{n+1}+a_{n+2}\label{101}.
\end{eqnarray}
\end{lemma}

\proof  The relation (\ref{001}) is a consequence of  the
consistence
$$
     p_{00^n 1}+p_{00^n0}=p_{00^n}
$$
and the facts $p_{00^n}=a_{n+1}$ and $p_{00^n0}=a_{n+2}$; the
relation (\ref{100}) is a consequence of  the invariance
$$
     p_{10^n 0}+p_{00^n0}=p_{0^n0}
$$
and the same facts; to obtain the relation (\ref{101}) we need
both the invariance and the consistence:
$$
  p_{10^n 1}+p_{00^n1}=p_{0^n1}=p_{0^n}-p_{0^n0}.
$$
Then by (\ref{001}) we get
$$
    p_{10^n 1}=p_{0^n}-p_{0^n0}-p_{00^n1} = a_n - 2a_{n+1} +a_{n+2}.
$$
$\Box$

Let $a, b, c$ be three positive numbers such that $a+b\ge c$.
Consider the function
\begin{equation}\label{F}
F(t, u, v, w)= h( t) + h(u) +h(v) + h(w)
\end{equation}
defined on $\mathbb{R}^{+4}$, where $h(x) = - x \log x$.

\begin{lemma}\label{maximize}
Under the condition
$$
   t+v=a,\quad u+w=b,\quad t+u=c
$$
the function $F$ defined by (\ref{F}) attains its maximum at
\begin{equation}\label{tuvw}
                               t= \frac{ ac}{a+b},\quad
u= \frac{ b c}{a+b},\quad v=\frac{ a (a + b-c)}{a+b}, \quad
w=\frac{b (a + b-c)}{a+b}.
\end{equation}
\end{lemma}
\proof From the condition we may solve $t, u, v$ as functions of
$w$:
$$
   t=c - b+w, \quad u=b-w, \quad v= a +b -c - w.
$$ So,  maximizing $F(t, u,v,w)$ under the condition
becomes maximizing the function
$$
F(w) = h(c - b+w)+h (b-w) +
        h(a +b -c - w) +h(w)
$$
which is strictly concave in its domain. Since $h'(x) = -1 - \log
x$, we have
$$
     F'(w) =- \log(c - b+w)+\log(b-w) +
        \log (a +b -c - w) -\log w.
$$
Solving $F'(w)=0$, we get $ w=\frac{b (a + b-c)}{a+b}. $ The
corresponding $t, u, v$ are as announced in (\ref{tuvw}) $\Box$



\begin{lemma}\label{unicite} Suppose that $\{a_n\}_{n\ge 0}$ is a
 sequence of real numbers such that
\begin{equation}
    1=a_0 \ge a_1 \ge a_2 \ge \dots; \quad a_i - 2 a_{i+1} +
    a_{i+2}\ge 0 \ \ (i\ge 0).
\end{equation}
There exists an invariant measure $\mu$ such that if $w$ is not a
block of $0$'s, we have
\begin{equation}\label{*w*}
p_{\epsilon w \epsilon'} = \frac{p_{\epsilon
w}p_{w\epsilon'}}{p_w} \qquad (\forall \epsilon, \epsilon' \in
\{0,1\}).
\end{equation}
The above recursion relations (\ref{*w*}) together with
(\ref{000}-\ref{101}) completely determine the measure $\mu$,
which is the unique maximal entropy measure among those invariant
measures $\nu$ such that $\nu([0^n]) = a_n$ for $n\ge 1$.
\end{lemma}

\proof For any $\mu in \mathcal M(\mathbf{a},\mathcal W)$, we must
have $\mu([0])=a_1$ and $\mu([1])=1-a_1$. By
Proposition~\ref{approximation}, we are led to find the measure
$\mu^{(n+2)} $ which maximizes $ h_{\mu}^{(n+2)}$ for each $n\ge 0
$. Let $\mu$ be an arbitrary invariant measure in
$\mathcal{M}(\mathbf{a}, \mathcal{W}_{n+2})$, $n\ge 0 $. We
identify $\mu$ with the sequence $p_w=\mu([w])$ indexed by finite
words. By Lemma
\ref{*0^k*}, we have \begin{eqnarray*} p_0  & = & a_1,\quad  p_1\ =1-a_1\\
p_{00} &=&a_2,\quad p_{01} = p_{10}= a_1-a_2,\qquad p_{11}= 1-2a_1
+a_2
\end{eqnarray*}
So, we have
\begin{equation}\label{entropy2}
     h_\mu^{(2)} = h(a_2) +2 h(a_1-a_2) - h(a_1) - h(1-a_1) + h(1-2 a_1 +
     a_2).
\end{equation}

Let us now consider the conditional entropy $h^{(n+2)}_\mu$ for
$n\ge 1$, which is  a function of $\mu$-measures of
$(n+2)$-cylinders, whose general form is $[\epsilon w \epsilon']$
with $w$ a word of length $n$ and $\epsilon, \epsilon' \in
\{0,1\}$. If $w=0^n$, by Lemma \ref{*0^k*}, the measures of the
four cylinders $[\epsilon 0^n\epsilon']$ with $\epsilon,
\epsilon'\in \{0, 1\}$ are determined by $\{a_k\}_{k\ge 1}$. If
$w\not=0^n$, the four quantities $p_{\epsilon w\epsilon'}$  are
linked to each other by
 $$
 p_{0w0}+p_{0w1}
     = p_{0w}, \quad
      p_{1w0}+p_{1w1} = p_{1w}, \quad
     p_{0w0}+p_{1w0}  = p_{w0}
$$
through measures of $(n+1)$-cylinders: $ a:=p_{w0}, b:=p_{1w},
c:=p_{0w}$. Consider the four measures $p_{\epsilon w\epsilon'}$
as variables, there is only one free variable and three others are
linked to it. Thus to any word $w\not=0^n$ of length $n$ is
associated a free variable.

In fact, $h_\mu^{(n+2)}$ is the sum of all terms
\begin{equation}\label{TT}
    -\sum_{\epsilon, \epsilon'} p_{\epsilon w \epsilon'}
    \log \frac{p_{\epsilon w\epsilon'}}{p_{\epsilon w}}
    \qquad (w \in\{0,1\}^n)
\end{equation}
If $w=0^n$, the corresponding term (\ref{TT}) is a constant
depending on the sequence $\{a_n\}$ (see Lemma~\ref{*0^k*}). If
$w\not=0^n$, there is a free variable in the term (\ref{TT}). So,
maximizing $h_\mu^{(n+2)}$ is equivalent to maximizing all above
terms, or equivalently to maximizing \begin{equation}\label{sum*w*
} -\sum_{\epsilon, \epsilon'} p_{\epsilon w \epsilon'} \log
p_{\epsilon w\epsilon'}.
\end{equation}
Applying Lemma \ref{maximize} to the above function in
(\ref{sum*w* }) provides us with the maximal point $\mu^{(n+2)} $
described by (\ref{*w*}) with $|w|=n $. It is easy to check that
the family $\{p_w\}$ defined by (\ref{000}-\ref{101}) and
(\ref{*w*}) for all words $w$, verifies the consistency and the
invariance conditions. The measure determined by the family is the
unique measure of maximal entropy and $\mu^{(n+2)} $ is its
$(n+1)$-Markov measure. $\Box$

\begin{lemma} The entropy of the invariant measure $\mu$ of maximal entropy determined in the
above lemma is equal to
$$
   h_\mu = -h(1-a_1) + \sum_{j=0}^\infty h(a_j - 2a_{j+1}+_{j+2})
$$
where $h(x) = -x \log x$.
\end{lemma}

\proof Recall that for any invariant measure $\mu$ we have
$$
     h_\mu =  \lim_{n\to \infty} h_\mu^{(n)} = \inf_n h_\mu^{(n)}
$$
where
$$
     h_\mu^{(n)}  = - \sum_{|w|=n-1} \ \ \sum_{\epsilon=0, 1}p_{w\epsilon} \log
     \frac{p_{w\epsilon}}{p_w}.
$$
Thus we may write
\begin{equation}\label{h-mu}
 h_\mu  =  h_\mu^{(2)} + \sum_{n=1}^\infty (h_\mu^{(n+2)} -
    h_\mu^{(n+1)})
\end{equation}

For $n \ge 0$,  write
\begin{eqnarray}
h_\mu^{(n+2)}
       & =&  - {\sum_{|w|=n}}' \sum_{\epsilon, \epsilon'} p_{\epsilon w\epsilon'} \log
     \frac{p_{\epsilon w\epsilon'} }{p_{\epsilon w}}
     - \sum_{\epsilon, \epsilon'} p_{\epsilon 0^n\epsilon'} \log
     \frac{p_{\epsilon0^n\epsilon'} }{p_{\epsilon0^n}} \nonumber \\
     &=&: I_1(n)
     +I_2(n)  \label{condition-entropy}
\end{eqnarray}
where $\sum'$ means that the sum is taken over $w\not=0^n$ ($0^0$
meaning the empty word so that $I_1(0) =0$). When $n=0$, we get
\begin{equation}\label{entropy2bis}
     h_\mu^{(2)} = h(a_2) +2 h(a_1-a_2) - h(a_1) - h(1-a_1) + h(1-2 a_1 +
     a_2).
\end{equation}
This coincides with (\ref{entropy2}).

Suppose $n\ge 1$. By the recursion relation (\ref{*w*}), we have
\begin{eqnarray}
   I_1(n)
      & = & - {\sum_{|w|=n}}' \sum_{\epsilon, \epsilon'} p_{\epsilon w\epsilon'} \log
     \frac{p_{w\epsilon'} }{p_{w}} \nonumber \\
     & = & -\sum_{|w|=n} \sum_{\epsilon, \epsilon'} p_{\epsilon w\epsilon'} \log
     \frac{p_{w\epsilon'} }{p_{w}}
     +  \sum_{\epsilon, \epsilon'} p_{\epsilon 0^n\epsilon'} \log
     \frac{p_{0^n\epsilon'} }{p_{0^n}} \nonumber \\
     & =& h_\mu^{(n+1)} + I_3(n) \label{I1}
\end{eqnarray}
where
$$
  I_3(n) =\sum_{\epsilon, \epsilon'} p_{\epsilon 0^n\epsilon'} \log
     \frac{p_{0^n\epsilon'} }{p_{0^n}}.
$$
From (\ref{condition-entropy}) and (\ref{I1}) we get
\begin{equation}\label{Entropy-difference}
h_\mu^{(n+2)} -h_\mu^{(n+1)} = I_2(n) + I_3(n).
\end{equation}
On the one hand, by using the invariance and the consistence we
can simplify $I_3$ as follows
\begin{eqnarray}\label{I3}
     I_3(n) & = & \sum_{ \epsilon'} p_{0^n\epsilon'} \log
     \frac{p_{0^n\epsilon'} }{p_{0^n}}\nonumber\\
     &=&
     \sum_{\epsilon'} p_{ 0^n\epsilon'} \log
     p_{0^n\epsilon'}  - p_{0^n} \log p_{0^n} \nonumber\\
      & =& a_{n+1}\log a_{n+1} + (a_n - a_{n+1}) \log (a_n - a_{n+1})
      - a_n \log a_n.\qquad
\end{eqnarray}
On the other hand, we have
\begin{eqnarray}\label{I2}
    I_2(n) & = &  - \sum_{\epsilon, \epsilon'} p_{\epsilon 0^n\epsilon'} \log
     p_{\epsilon0^n\epsilon'}  + \sum_{\epsilon} p_{\epsilon0^n}
     \log p_{\epsilon0^n}\nonumber \\
     & = & - a_{n+2} \log a_{n+2} - 2 (a_{n+1}-a_{n+2}) \log
     (a_{n+1}-a_{n+2}) \nonumber \\
     & & \quad - (a_n -2 a_{n+1}+a_{n+2})\log (a_n -2
     a_{n+1}+a_{n+2}) \nonumber\\
     &  &   \quad + a_{n+1} \log a_{n+1} + (a_n -a_{n+1}) \log (a_n
     -a_{n+1}).
\end{eqnarray}

By combining (\ref{Entropy-difference}), (\ref{I3}) and (\ref{I2})
we
get \begin{eqnarray*} \varphi(n) &:= & h_\mu^{(n+2)} -h_\mu^{(n+1)}\\
&=& [h(a_{n+2}) - 2 h(a_{n+1}) + h(a_n)] + 2 [h(a_{n+1}-a_{n+2})
-h(a_n - a_{n+1})] \\ & & \ \ \ \ \ + h(a_n - 2a_{n+1}+a_{n+2}).
\end{eqnarray*}
Finally using (\ref{h-mu}) we  get
\begin{eqnarray*}
    h_\mu
     =  h_\mu^{(2)} + \sum_{n=1}^\infty \varphi(n)
         =  -h(1-a_1) + \sum_{j=0}^\infty h(a_j -
        2a_{j+1}+a_{j+2}).
\end{eqnarray*}
$\Box$
\medskip

Remark that the measure of maximal entropy in $R(a, \mathcal{W})$
is not necessarily ergodic. Here is an example. If $a_n =a$
($\forall n \ge 1$) is constant, then there is a unique invariant
measure in $\mathcal{M}(a, \mathcal{W})$, which is $a
\delta_{\bar{0}} + (1-a) \delta_{\bar{1}}$.  In this case, $R(a,
\mathcal{W})$ is not empty but of zero  entropy.  Notice that
$R(a, \mathcal{W})$ contains  no point in the support of the
unique invariant measure. If $a=\frac{1}{2}$, $R(\frac12,
\mathcal{W})$ contains the following point
$$
   01001100011100001111...
$$
(The terms in the two sequences $\{0^k\}$ and $\{1^k\} $ are
alternatively appended.)

We thank the referee for telling us that Pfister and Sullivan
\cite{PfSu} also obtained the same result as our Theorem
\ref{specification-saturated} but with a different method.

\end{document}